\documentclass[reqno,a4paper,11pt]{amsart}

\usepackage{amsmath}
\usepackage{amssymb}
\usepackage{amsthm}
\usepackage{amscd}
\usepackage{mathrsfs}

\theoremstyle{plain}
\newtheorem{thm}{Theorem}[section]
\newtheorem{lem}[thm]{Lemma}
\newtheorem{prop}[thm]{Proposition}
\newtheorem{cor}[thm]{Corollary}
\theoremstyle{definition}
\newtheorem{defn}[thm]{Definition}

\newtheorem{rem}[thm]{Remark}

\newcommand{\ri}{\mathfrak{o}}
\newcommand{\mi}{\mathfrak{p}}

\newcommand{\rad}{\mathfrak{a}}
\newcommand{\seqrad}{\mathfrak{a}}

\newcommand{\g}{\mathfrak{g}}
\newcommand{\Z}{\mathbf{Z}}
\newcommand{\C}{\mathbf{C}}
\newcommand{\He}{\mathcal{H}}

\newcommand{\Irr}{\mathrm{Irr}}

\newcommand{\Lamdba}{\Lambda}

\newcommand{\ad}{\mathrm{ad}}
\newcommand{\Ad}{\mathrm{Ad}}

\newcommand{\e}{\sqrt{\varepsilon}}
\newcommand{\p}{\varpi}

\newcommand{\J}{\mathfrak{J}}

\title{Representations of unramified $U(2,2)$ over 
a $p$-adic field I:  representations
of non-integral level}
\author{Michitaka Miyauchi}
\address{
School of Mathematics,
University of East Anglia,
Norwich NR4 7TJ UK
}
\email{M.Miyauchi@uea.ac.uk}
\begin{document}

\maketitle

\pagestyle{myheadings}
\markboth{MICHITAKA MIYAUCHI}{REPRESENTATIONS OF $U(2,2)$ OVER
A $p$-ADIC FIELD I}

\begin{abstract}
Let $F_0$ be a non-archimedean local field of 
odd residual characteristic
and let $G$ be the unramified unitary group $U(2,2)$
defined over $F_0$.
In this paper,
we give a classification
 of the irreducible smooth representations of $G$
of non-integral level
using the Hecke algebraic 
method developed by Allen Moy for $GSp(4)$.
\end{abstract}

\section*{Introduction}
Let $F_0$ be a non-archimedean local field
of odd residual characteristic
and let $G$ be the unramified unitary group
$U(2,2)$ defined over $F_0$.
Although Konno \cite{Konno} classified the non-supercuspidal 
representations of $G$,
supercuspidal representations of $G$ has not been classified.
The purpose of this paper is to 
classify the
irreducible smooth representations of $G$ of non-integral level.

In \cite{U21} and \cite{GSp4},
Moy gave a classification of 
the irreducible smooth representations 
of unramified $U(2,1)$ and $GSp(4)$ over $F_0$,
based on the concepts of
nondegenerate representations
and Hecke algebra isomorphisms.
A nondegenerate representation of $GSp(4)$
is an irreducible representation $\sigma$ 
of an open compact subgroup $K$
which satisfies a certain cuspidality or semisimplicity condition.
An important property of nondegenerate representations 
of $GSp(4)$ is that
every irreducible smooth representation of $GSp(4)$
contains some nondegenerate representation.
For $\sigma$ a nondegenerate representation of
$GSp(4)$,
the set of equivalence classes of irreducible representations
of $GSp(4)$
which contain $\sigma$
can be identified with
the set of equivalence classes of irreducible representations of a Hecke algebra $\He$ associated
to $\sigma$.
Moy described $\He$ as a Hecke algebra of
some smaller group
and thus
reduced the classification of 
the irreducible smooth representations of $GSp(4)$ which contain
$\sigma$
to that of a smaller group.

In this paper, we  attempt to  classify
the irreducible smooth representations of $G$ 
by this method.
The keypoint of our classification is 
to construct an analog of nondegenerate representations
of $GSp(4)$ for $G$.

In \cite{MP},
Moy and Prasad developed
the concept of nondegenerate representations 
into that of unrefined minimal $K$-types for reductive 
$p$-adic groups.
For classical groups,
Stevens \cite{St3} gave an explicit construction of 
unrefined minimal $K$-types as fundamental skew strata,
based on the results of Bushnell and Kutzko \cite{BK2}
and Morris \cite{Morris-1}.
However,
neither of these give a classification of 
the irreducible smooth representations of 
$p$-adic classical groups.

Throughout this paper,
we use the notion of fundamental skew strata introduced by
\cite{St3} for our nondegenerate representations of $G$.
Let $F$ be the unramified quadratic extension over $F_0$.
Then $G$ is realized as the group of isometries
of an $F/F_0$-hermitian form on 4-dimensional 
$F$-vector space $V$.
According to \cite{St3},
a skew stratum 
is a 4-tuple $[\Lambda, n, r, \beta]$.
A periodic lattice function $\Lambda$ with a certain duality
induces
a filtration $\{P_{\Lambda, k}\}_{k \geq 1}$ on 
a parahoric subgroup $P_{\Lambda, 0}$ of $G$.
Integers $n > r \geq 0$ 
and an element $\beta$ in the Lie algebra of $G$
determine a character $\psi$ of the group $P_{\Lambda, r+1}$
which is trivial on $P_{\Lambda, n+1}$.
Writing $e(\Lambda)$ for the period of $\Lambda$,
we refer to $n/e(\Lambda)$ as the level of the stratum.

In Section~\ref{strata},
we prove Theorem~\ref{thm:strict},
a rigid result on the existence of fundamental skew strata.
This theorem,
which is an analog of the result in \cite{GSp4},
says that
every irreducible representation of $G$
contains some fundamental skew stratum $[\Lambda, n, n-1, \beta]$ such that $\{ P_{\Lambda, n}\}_{n \geq 0}$ is 
the standard filtration of $P_{\Lambda, 0}$.
So we can start our classification with 
7 filtrations up to conjugacy.

In the latter part of this paper,
we give a classification of the irreducible 
representations of $G$ which contain a
skew stratum $[\Lambda, n, n-1, \beta]$ 
whose level is not an integer.
As in \cite{St3},
we may further assume that $\beta$ is a semisimple element in 
$\mathrm{Lie}(G)$ by
replacing $P_{\Lambda, n}$ with a nonstandard filtration subgroup
of a parahoric subgroup.

In Section~\ref{supercuspidal},
we consider the case when the $G$-centralizer of $\beta$
forms a maximal compact torus in $G$.
Stevens \cite{St2} and \cite{St1}
gave a method to construct irreducible supercuspidal representations of $p$-adic classical groups from such a stratum.
Applying his result, we can construct the irreducible representations
of $G$ containing such a stratum.
Those are irreducibly induced from open and compact subgroups.

In Section~\ref{Hecke},
we treat the case when the $G$-centralizer $G'$ of $\beta$
is not compact.
We describe the Hecke algebra associated to such a stratum 
as 
a certain Hecke algebra of $G'$.
This allows us to identify the equivalence classes of 
irreducible representations of $G$ containing such a stratum
with those of irreducible representations of $G'$.
The construction of Hecke algebra isomorphisms here
is along the lines of that in \cite{U21} and \cite{GSp4}.
This method is valid only if $G'$ is tamely ramified over $F_0$.
But there are no wildly ramified $G'$
for our group $G =U(2, 2)$.
This is the reason why Moy's classification works well for $G$.

Moreover, we consider the condition when
two fundamental skew strata occur in a common 
irreducible smooth representation in $G$.
This problem relates to  the
\lq\lq intertwining implies conjugacy" property.
We prove that this
property holds 
among  the skew strata $[\Lambda, n, n-1, \beta]$
with compact centralizers.
In the other cases,
we normalize skew strata
to be appropriate to this problem.
Then the
\lq\lq intertwining implies conjugacy" property 
of normalized skew strata is trivial
except only one case (case (\ref{sec:h_int}d)).

The remaining problem that we need to consider is to classify the irreducible 
representations of $G$ of integral level.
The level 0 representations have been classified by
Moy-Prasad \cite{Jacquet} and Morris \cite{Morris0} for more 
general $p$-adic reductive groups.
In a second part of this article,
we will complete  our classification of 
the irreducible representations of $G$,
by classifying those
containing a fundamental skew stratum of positive integral level.

The Hecke algebra isomorphisms established for $G$
preserve the unitary structure of Hecke algebras,
so we can calculate
formal degrees
of supercuspidal representations of $G$ via those isomorphisms.
We hope to return to this in the future.

Part of this article is based on the author's 
doctoral thesis.
The author would like to thank his supervisors, Tetsuya Takahashi 
and Tadashi Yamazaki, for 
their useful comments and 
patient encouragement during this work.
The author also thank Kazutoshi Kariyama and Shaun Stevens
for corrections and helpful comments.
The research for this paper was partially supported by EPSRC grant
GR/T21714/01.

\section{Preliminaries}\setcounter{equation}{0}
In this section, we recall the notion of fundamental skew strata
for unramified unitary groups over a non-archimedean
local field.
For details in more general settings,
one should consult \cite{BK1}, \cite{BK2} and \cite{St2}.

\subsection{Filtrations}
Let $F_0$ be a non-archimedean local field
of odd residual characteristic.
Let $\ri_0$ denote the ring of integers in $F_0$,
$\mi_0$ the prime ideal in $\ri_0$,
$k_0 = \ri_0/\mi_0$ the residue field,
and $q$ the number of elements in $k_0$.

Let $F = F_0[\e]$, 
$\varepsilon$ a nonsquare element in $\ri_0^\times$,
denote the unramified quadratic extension over $F_0$.
We denote by
$\ri_F$, $\mi_F$, $k_F$ the objects for $F$ 
analogous to those above for $F_0$.
For $x$ in $F$ or $k_F$,
we write $\overline{x}$ for the Galois conjugate of $x$.
Since $F$ is unramified over $F_0$,
we can (and do) select a common uniformizer $\p$ 
of $F_0$ and $F$.

Let $V$ be an $N$-dimensional $F$-vector space 
equipped with a nondegenerate hermitian form $f$
with respect to $F/F_0$.
We put
$A = \mathrm{End}_F(V)$ and 
$\widetilde{G} = A^\times$.
Let $\sigma$ denote
the involution on $A$ induced by $f$.
We also put
$G= \{ g \in \widetilde{G}\ |\ g\sigma(g) = 1 \}$,
the corresponding unramified unitary group over $F_0$,
and
$\g = \{ X \in A\ |\ X+\sigma(X)=0 \}
\simeq \mathrm{Lie}(G)$.

Recall from \cite{BK2} (2.1) that
an $\ri_F$-lattice sequence in $V$ is a function $\Lambda$ from $\mathbf{Z}$
to the set of $\ri_F$-lattices in $V$ such that
\begin{enumerate}
\item[(i)] $\Lambda(i) \supset \Lambda(i+1)$, $i \in \mathbf{Z}$;

\item[(ii)]
there exists an integer $e(\Lambda)$ 
called the $\ri_F$-period of $\Lambda$
such that
$\varpi \Lambda(i) = \Lambda(i+e(\Lambda))$, $i \in \mathbf{Z}$.
\end{enumerate}
We say that
an $\ri_F$-lattice sequence $\Lambda$ is strict 
if $\Lambda$ is injective.

For $L$ an $\ri_F$-lattice in $V$,
we define its dual lattice $L^\#$ by
$L^\# = \{ v \in V\ |\ f(v, L) \subset \ri_F\}$.
An $\ri_F$-lattice sequence $\Lambda$ in $V$ is called
self-dual if there exists an integer $d(\Lambda)$
such that
$\Lambda(i)^\# = \Lambda(d(\Lambda)-i)$,
$i \in \Z$. 

An $\ri_F$-lattice sequence $\Lambda$ in $V$
induces
a filtration $\{\rad_n(\Lambda)\}_{n \in \Z}$ on $A$ by
\[
\mathfrak{a}_n(\Lambda)
= \{ X \in A\ |\ X\Lambda(i) \subset \Lambda(i+n),\ i \in 
\mathbf{Z}\},\ n \in \Z.
\]
This filtration 
determines a sort of \lq\lq valuation" map
$\nu_\Lambda$ on $A$ by
\[
\nu_\Lambda(x) = \sup\{ n \in \Z\ |\ x \in \seqrad_n(\Lambda)\},\ 
x \in A\backslash \{0\},
\]
with the usual understanding that
$\nu_\Lambda(0) = \infty$.

For $\Lambda$ an $\ri_F$-lattice sequence in $V$
and $k \in \Z$,
we define a translate $\Lambda +k$ of $\Lambda$ by
$(\Lambda +k)(i) = \Lambda(i+k)$, $i \in \Z$.
Then we have $\seqrad_n(\Lambda) = \seqrad_{n}(\Lambda +k)$, $n \in \Z$.
For $g \in G$
and  
$\Lambda$ a self-dual $\ri_F$-lattice sequence in $V$,
we define a self-dual $\ri_F$-lattice sequence $g\Lambda$
by $(g\Lambda)(i) = g\Lambda(i)$, $i \in \Z$.
Note that
$e(g\Lambda) = e(\Lambda)$ and $d(g\Lambda) = 
d(\Lambda)$.
\begin{rem}
Given an $\ri_F$-lattice sequence $\Lambda$ in $V$,
we define its dual sequence $\Lambda^\#$ by
$\Lambda^\#(i) = \Lambda(-i)^\#$, $i \in \Z$.
Since 
$\seqrad_n(\Lambda^\#) = \sigma(\seqrad_n(\Lambda))$, $n \in \Z$,
it follows from \cite{BK2} (2.5) that
$\Lambda$ is self-dual if and only if
$\sigma(\seqrad_n(\Lambda)) = \seqrad_n(\Lambda)$,
$n \in \Z$. 
\end{rem}
For $\Gamma$ an $\ri_F$-lattice in $A$,
we define its dual by
$\Gamma^*  =  \{ X \in A\ |\ \mathrm{tr}_{A/F_0}(X\Gamma) \subset \mi_0\}$,
where 
$\mathrm{tr}_{A/F_0}$ denotes the composition of traces
$\mathrm{tr}_{F/F_0}\circ \mathrm{tr}_{A/F}$.
By \cite{BK2} (2.10),
if $\Lambda$ is an $\ri_F$-lattice sequence in $V$,
then we have 
$\seqrad_n(\Lambda)^* = \seqrad_{1-n}(\Lambda)$.

For $\Lambda$ a self-dual $\ri_F$-lattice sequence in $V$,
 we will write
$\g_{\Lambda, n} = \g \cap \seqrad_{n}(\Lambda)$,
$n \in \Z$.
We also set
$P_{\Lambda, 0} = 
G\cap \seqrad_0(\Lambda)$
and 
$P_{\Lambda, n} = G\cap (1 + \seqrad_{n}(\Lambda))$,
$n \geq 1$.
Then $\{ P_{\Lambda, n}\}_{n \geq 0}$
is a filtration of the parahoric subgroup $P_{\Lambda, 0}$
of ${G}$ by its open normal subgroups.

We fix an additive character
$\Omega$ of $F_0$ with conductor $\mi_0$.
Let 
${}^\wedge$ denote the Pontrjagin dual
and for
$x$ a real number, 
let $[x]$ denote the greatest integer less than or equal to $x$.
\begin{prop}\label{thm:Morris}
\label{prop:p}
Let $\Lambda$ be a self-dual $\ri_F$-lattice sequence in  
$V$ and 
let 
$n, r \in \Z$ satisfy
$n > r \geq [n/2] \geq 0$.
Then
the map $p \mapsto p-1$ induces an
isomorphism
$P_{\Lambda, r+1}/P_{\Lambda, n+1} \simeq 
\g_{\Lambda, r+1}/\g_{\Lambda, n+1}$ and
there exists an isomorphism 
\[
\g_{\Lambda, -n}/\g_{\Lambda, -r} \simeq
(P_{\Lambda, r+1}/P_{\Lambda, n+1})^\wedge;
b+ \g_{\Lambda, -r} \mapsto \psi_b,
\]
where
$\psi_b(p) =
\Omega(\mathrm{tr}_{A/F_0}(b(p-1)))$, 
$p \in P_{\Lambda, r+1}$.
\end{prop}
\subsection{Skew strata}
\begin{defn}[\cite{BK2} (3.1), \cite{St2} Definition 4.5]
(i) A stratum in $A$ is a 4-tuple $[\Lambda, n, r, \beta]$
consisting of an $\ri_F$-lattice sequence $\Lambda$ in $V$,
integers $n, r $ verifying $n > r \geq 0$, 
and an element $\beta$ in $\seqrad_{-n}(\Lambda)$.
We say that two strata $[\Lambda, n, r, \beta_i]$, $i = 1,2$,
are equivalent if $\beta_1 - \beta_2 
\in \seqrad_{-r}(\Lambda)$.

(ii) A stratum $[\Lambda, n, r, \beta]$ in $A$ is called skew
if $\Lambda$ is self-dual and $\beta \in \g_{\Lambda, -n}$.
\end{defn}

The fraction 
$n/e(\Lambda)$ is called the level of the stratum.
If $n > r \geq [n/2]$,
then by
Proposition~\ref{thm:Morris},
the equivalence classes of skew strata 
of the form 
$[\Lambda, n, r, \beta]$ parametrize the dual of
$P_{\Lambda, r+1}/P_{\Lambda, n+1}$.

For
$[\Lambda, n, r, \beta]$ a stratum in $A$,
we set
$y_\beta = \p^{n/k} \beta^{e(\Lambda)/k}$,
where
$k = (e(\Lambda), n)$.
Then $y_\beta$ belongs to $\rad_0(\Lambda)$
and 
its characteristic polynomial 
$\Phi_\beta(X)$
lies in $\ri_F[X]$.
We define the characteristic polynomial 
$\phi_\beta(X)$ of this stratum 
to be the reduction  modulo $\mi_F$ of $\Phi_\beta(X)$.
\begin{defn}[\cite{BK1} (2.3)]
A
stratum $[\Lambda, n, r, \beta]$ in $A$ is called fundamental 
if
$\phi_\beta(X) \neq X^N$.
\end{defn}

The representations of $G$ we will consider are
always assumed to be smooth and complex.
Let 
$\pi$ be a smooth representation of $G$
and $[\Lambda, n, r, \beta]$
a skew stratum with $n > r \geq [n/2]$.
We say that $\pi$ contains 
$[\Lambda, n, r, \beta]$ if 
the restriction of $\pi$ to 
$P_{\Lambda, r+1}$ contains the character
$\psi_\beta$.
A smooth representation $\pi$ of $G$
is called of positive level
if $\pi$ has no non-zero $P_{\Lambda, 1}$-fixed vector,
for any self-dual $\ri_F$-lattice sequence $\Lambda$ in $V$.
By \cite{St3} Theorem 2.11,
an irreducible smooth representation of $G$ of positive level
contains a fundamental skew stratum $[\Lambda, n, n-1, \beta]$.

For $g \in \widetilde{G}$ and $x \in A$,
we write $\Ad(g)(x) = gxg^{-1}$.
\begin{prop}\label{thm:MP}\label{prop:MP}
Let $[\Lambda, n, r, \beta]$ and $[\Lambda', m, s, \gamma]$
be skew strata in $A$
contained in some irreducible smooth representation of $G$.

(i) There exists $g \in G$ such that
$(\beta +\g_{\Lambda, -r})\cap \Ad(g)(\gamma +\g_{\Lambda', -s}) \neq \emptyset$.

(ii)  If $[\Lambda, n, r, \beta]$ is fundamental,
then we have $n/e(\Lambda) \leq m/e(\Lambda')$.

(iii) If 
$n/e(\Lambda) =m/e(\Lambda')$,
then $\phi_\beta(X) = \phi_\gamma(X)$.
\end{prop}
\begin{proof}
The proof is very similar to those of 
\cite{HM2} Theorem 4.1 and Corollary 4.2.
\end{proof}
By Proposition~\ref{prop:MP},
we can define the level and the characteristic polynomial of 
an irreducible smooth representation $\pi$ of $G$ of positive
level
to be those of the fundamental skew strata contained in $\pi$.

Given a skew stratum $[\Lambda, n, r, \beta]$ in $A$,
we define its formal intertwining
\[
I_{G}[\Lambda, n, r, \beta]
= \{ g \in G\ |\ (\beta + 
\g_{\Lambda, -r})
\cap \Ad(g)(\beta + \g_{\Lambda, -r}) \neq \emptyset \}.
\]
If $n > r \geq [n/2]$, then it is nothing other than
the intertwining of the character $\psi_\beta$ of $P_{\Lambda, r+1}$
in $G$.

\subsection{Semisimple strata}
\begin{defn}[\cite{BK1} (1.5.5), \cite{BK2} (5.1)]
A stratum 
$[\Lambda, n, r, \beta]$ in $A$ is called simple if

(i) the algebra $E = F[\beta]$ is a field,
and $\Lambda$ is an $\ri_E$-lattice sequence;

(ii) $\nu_\Lambda(\beta) = -n$;

(iii) $\beta$ is minimal over $F$
in the sense of \cite{BK1} (1.4.14).
\end{defn}
\begin{rem}
The definition above is restricted for our use.
With the notion of \cite{BK1},
our simple stratum 
is a simple stratum $[\Lambda, n, r, \beta]$
with $k_0(\beta, \Lambda) = -n$.
\end{rem}

The following criterion of the simplicity of strata
is well known.
\begin{prop}[\cite{char_l} Proposition 1.5]\label{prop:pure}
Let
$\Lambda$ be a strict $\ri_F$-lattice sequence in $V$ 
of period $N$
and $n$ an integer coprime to  $N$.
Suppose that
$[\Lambda, n, n-1, \beta]$ is a fundamental stratum in $A$.
Then $F[\beta]$ is a totally ramified extension of degree $N$
over $F$ and
$[\Lambda, n, r, \beta]$ is simple,
for any $n > r \geq 0$.
\end{prop}

Let $[\Lambda, n, r, \beta]$ be a stratum in $A$.
We assume that there is a non-trivial $F$-splitting $V = V^1 \oplus V^2$ such that
\begin{enumerate}
\item[(i)]
$\Lambda(i) = \Lambda^1(i) \oplus \Lambda^2(i)$,
$i \in \Z$,
where $\Lambda^j(i) = \Lambda(i) \cap V^j$, for $j = 1, 2$;

\item[(ii)] $\beta V^j \subset V^j$, $j =1,2$.
\end{enumerate}
For $j = 1,2$, we write $\beta_j = \beta|_{V^j}$.
By \cite{BK2} (2.9),  we get
a stratum $[\Lambda^j, n, r, \beta_j]$
 in $\mathrm{End}_F(V^j)$.
Recall from \cite{BK2} (3.6) that
a stratum $[\Lambda, n, r, \beta]$ in $A$ is called split if
\begin{enumerate}
\item[(iii)] $\nu_{\Lambda^1}(\beta_1) = -n$ and
$X$ does not divide $\phi_{\beta_1}(X)$;

\item[(iv)] either $\nu_{\Lambda^2}(\beta_2) > -n$, or else all the following 
conditions hold:

\begin{enumerate}
\item[(a)] $\nu_{\Lambda^2}(\beta_2) = -n$ and
$X$ does not divide $\phi_{\beta_2}(X)$,

\item[(b)] $(\phi_{\beta_1}(X), \phi_{\beta_2}(X)) = 1$.
\end{enumerate}
\end{enumerate}
\begin{defn}[\cite{St2} Definition 4.8, \cite{St3} Definition 2.10]
(i) (Inductive definition on the dimension of $V$)
 A stratum $[\Lambda, n, r, \beta]$ is called semisimple
if it is simple, or else
it is split as above
and satisfies the following conditions:

(a) $[\Lambda^1, n, r, \beta_1]$ is simple;

(b) 
$[\Lambda^2, n_2, r, \beta_2]$ is semisimple 
or $\beta_2 = 0$, where
$n_2 = \mathrm{max}\{-\nu_{\Lambda^2}(\beta_2),
r+1\}$.

\noindent
(ii)
A skew stratum in $A$ is called split
if it is split with respect to an orthogonal $F$-splitting
$V = V^1 \bot V^2$.
\end{defn}

If a skew stratum $[\Lambda, n, r, \beta]$ is split
with respect to $V = V^1 \bot V^2$,
then 
$\Lambda^j$ is a self-dual $\ri_F$-lattice sequence
in $(V^j, f|_{V^j})$ with $d(\Lambda^j) = d(\Lambda)$,
for $j =1, 2$.
\subsection{Hecke algebras}
Let $G$ be a unimodular, 
locally compact, totally disconnected topological group.
Let $J$ be an open compact subgroup of $G$
and let $(\sigma, W)$ be an irreducible smooth representation of $J$.
For $g \in G$, we write
$\sigma^g$ for the representation of
$J^g = g^{-1}Jg$ defined by
$\sigma^g (x) = \sigma(g x g^{-1})$, $x \in J^g$.
We define the intertwining of $\sigma$ in $G$ by
\[
I_G(\sigma) = 
\{ g \in {G}\ 
|\ \mathrm{Hom}_{J\cap J^g}(\sigma, \sigma^g) \neq 0 \}.
\]

Let $(\widetilde{\sigma}, \widetilde{W})$ denote 
the contragradient representation of $(\sigma, W)$.
The Hecke algebra $\He(G//J, \sigma)$ 
is the set of 
compactly supported functions $f : G  \rightarrow \mathrm{End}_{\C}(\widetilde{W})$
such that
\[
f(k g k') = \widetilde{\sigma}(k)f(g) \widetilde{\sigma}(k'),\
k, k' \in J,\ g \in G.
\]
Let $dg$ denote the Haar measure on $G$ normalized
so that the volume $\mathrm{vol}(J)$ of $J$ is $1$.
Then $\He(G//J, \sigma)$ becomes an algebra
under convolution relative to $dg$.
Recall from \cite{BK1} (4.1.1) that the support 
of $\He(G//J, \sigma)$
is the intertwining of  $\widetilde{\sigma}$ in $G$,
that is, 
\[
I_G(\widetilde{\sigma}) = 
\bigcup_{f\in \He(G//J, \sigma)}\mathrm{supp}(f).
\]

Since $J$ is compact,
there exists a $J$-invariant, positive definite hermitian form 
on $\widetilde{W}$.
This form induces an involution $X \mapsto \overline{X}$
on $\mathrm{End}_{\C}(\widetilde{W})$.
For $f \in \He(G//J, \sigma)$,
we define $f^* \in \He(G//J, \sigma)$ by 
$f^*(g) = \overline{f(g^{-1})},\ g \in G$.
Then the map $*: \He(G//J, \sigma) \rightarrow 
\He(G//J, \sigma)$ is an involution on $\He(G//J, \sigma)$.

Let 
$\Irr(G)$ denote the set of equivalence classes of 
irreducible smooth representations of $G$
and $\Irr(G)^{(J, \sigma)}$
the subset of $\Irr(G)$ consisting of
the elements
those $\sigma$-isotypic components are not zero.
Let $\Irr \He(G//J, \sigma)$ denote the 
set of equivalence classes of irreducible representations
of $\He(G//J, \sigma)$.
Then, by \cite{BK1} (4.2.5), there is a bijection 
$\Irr(G)^{(J, \sigma)} \simeq \Irr \He(G//J, \sigma)$.

\subsection{Unramified hermitian forms
on a 2-dimensional space}
We shall consider 
hermitian forms on the space
$V = F^2$.
Let $\{e_1, e_2\}$ be an $F$-basis of $V$.
Up to isometries,
there are two $F/F_0$-hermitian forms $f_0$
and $f_1$ on $V$, where
\begin{eqnarray*}
f_0(x, y) = 
\overline{x}_{1}y_1 +
\p \overline{x}_{2}y_2,\
f_1(x, y) = 
\overline{x}_{1}y_2 +
\overline{x}_{2}y_1, 
\end{eqnarray*}
for $x = x_1 e_1 + x_2 e_2$, $y = y_1 e_1 + y_2 e_2$
in $V$.

The space $(V, f_0)$ is anisotropic and
has the unique $\ri_F$-lattice $L$ such that $L^\# 
\supset L \supset \p L^\#$.
In fact, we have
$L = \ri_F e_1 \oplus \ri_F e_2$.
Therefore, every strict self-dual $\ri_F$-lattice sequence in $(V, f_0)$
is a translate of the following sequence $\Lambda$:
\begin{eqnarray}
\Lambda(2i) = \p^i L,\ 
\Lambda(2i+1) = \p^{i+1} L^\#,\ i \in \Z.
\end{eqnarray}
Hence
every strict self-dual $\ri_F$-lattice sequence $\Lambda'$
in $(V, f_0)$ satisfies
$e(\Lambda') = 2$ and $d(\Lambda')$ is odd.

Next, we consider the space $(V, f_1)$.
Its anisotropic part is trivial.
Set
\begin{eqnarray*}
N_0 = \ri_F e_1 \oplus \ri_F e_2,\
N_1 = \ri_F e_1 \oplus \mi_F e_2.
\end{eqnarray*}
Let $\Lambda$ be a strict self-dual $\ri_F$-lattice sequence
in $(V, f_1)$ of $\ri_F$-period 2.
Then there exist $g \in G$ and $k \in \Z$
such that
\begin{eqnarray}
(g\Lambda+k)(2i) = \p^i N_0,\ (g\Lambda+k)(2i+1) = 
\p^i N_1,\
i \in \Z.
\end{eqnarray}
In particular, the integer $d(\Lambda)$ should be even.

As a consequence, we obtain the following:
\begin{lem}\label{lem:anisotropic}
Let $f$ be an unramified hermitian form on a 
2-dimensional $F$-space $V$
and let $\Lambda$ be 
a strict self-dual $\ri_F$-lattice sequence $\Lambda$ in $V$ 
with $e(\Lambda) = 2$.
Then the parity of $d(\Lambda)$
is independent of the choice of $\Lambda$.
Moreover,
The space $(V, f)$ is anisotropic
if and only if
$d(\Lambda)$ is odd.
\end{lem}

\section{Fundamental strata for $U(2,2)$}
\setcounter{equation}{0}\label{strata}
In this section, 
we realize the unramified unitary group $U(2, 2)$ 
and state a rigid version of 
the existence of  fundamental skew strata
for $U(2,2)$, which implies that
every irreducible smooth representation of $U(2,2)$ of positive level
contains a fundamental skew stratum 
$[\Lambda, n, n-1, \beta]$ such that
$\Lambda$ is a strict $\ri_F$-lattice sequence.
\subsection{Reduction to strict lattice sequences}
From now on,
we will denote by $V$ 
the four dimensional space of column vectors $F^4$.
We write
$A = M_4(F)$ and
$\widetilde{G} = GL_4(F)$.
Let $e_i$ ($1 \leq i \leq 4$) be the standard basis vectors
of $V$,
and let $E_{i j}$ denote the element in $A$ whose $(k, l)$
entry is $\delta_{ik}\delta_{jl}$.
Set
\[
H = E_{14} + E_{23} + E_{32} + E_{41},
\]
and define a nondegenerate hermitian form $f$ on $V$
by
$f(x, y) = {}^t \overline{x} H  y$, $x, y \in V$.
The form $f$ induces an involution 
$\sigma$ on $A$  by
$\sigma(X) = H^{-1} {}^t\overline{X}H$, $X \in A$.

We define
$G = U(2,2) = \{ g \in \widetilde{G}\ |\ g\sigma(g) = 1 \}$
and
$\g = \{ X \in A\ |\ X +\sigma(X) = 0\} \simeq 
\mathrm{Lie}(G)$.
Then $\g$ consists of matrices of the form
\[
 \left(
\begin{array}{cccc}
Y & Z & C & a\sqrt{\varepsilon}\\
M & N & b\sqrt{\varepsilon} & -\overline{C} \\
D & c\sqrt{\varepsilon} & -\overline{N} & -\overline{Z} \\
d\sqrt{\varepsilon} & -\overline{D} & -\overline{M} & -\overline{Y}
\end{array}
\right),\ 
C, D, M, N, Y, Z \in F,\ a,b,c,d \in F_0.
\]

We recall the structure of strict self-dual $\ri_F$-lattice
sequences in $(V, f)$ from \cite{Morris-2} \S 1.
We define $\ri_F$-lattices $N_0$, $N_1$, and $N_2$ in $V$ by
\begin{eqnarray*}
& N_0 = \ri_F e_{1} \oplus \ri_F e_{2} \oplus \ri_F e_3 \oplus \ri_F e_4,\\
& N_1 = \ri_F e_{1} \oplus \ri_F e_{2} \oplus \ri_F e_3 \oplus \mi_F e_4,\\
& N_2 = \ri_F e_{1} \oplus \ri_F e_{2} \oplus \mi_F e_3 \oplus \mi_F e_4.
\end{eqnarray*}
Then we obtain a sequence of lattices
\begin{eqnarray*}
N_0^\# = N_0 \supsetneq N_1 \supsetneq N_2 = \p N_2^\# \supsetneq
\p N_1^\# \supsetneq \p N_0.
\end{eqnarray*}
A self-dual $\ri_F$-lattice sequence $\Lambda$ in $V$
is called standard if its image
$\Lambda(\Z) = \{\Lambda(i)\ |\ i \in \Z\}$ 
is 
contained in the set
$\{ \p^m N_0,\ \p^m N_1,\ \p^m N_2,\ \p^m N_1^\#\ |\ m \in \Z\}$.
By \cite{Morris-2} Proposition 1.10,
every self-dual $\ri_F$-lattice sequence is a $G$-conjugate of 
some standard sequence.

Up to translation, 
the standard strict self-dual $\ri_F$-lattice sequences
correspond
to the non-empty subsets $S $ of $\{ N_0, N_1, N_2\}$.
Therefore,
there are just following 7 standard strict self-dual $\ri_F$-lattice sequences in $V$ up to translation:

(1) $S = \{N_0, N_1, N_2\}$:
\begin{eqnarray*}
\Lambda(4i) = \p^i N_0,\
\Lambda(4i+1) = \p^i N_1,
\end{eqnarray*}
\begin{eqnarray}\label{eq:st_4}
\Lambda(4i+2) = \p^i N_2,\ 
\Lambda(4i+3) = \p^{i+1} N_1^\#,\ i \in \Z,
\end{eqnarray}

(2) $S = \{N_0, N_1\}$:
\begin{eqnarray}\label{eq:st_30}
\Lambda(3i) = \p^i N_0,\
\Lambda(3i+1) = \p^i N_1,\ 
\Lambda(3i+2) = \p^{i+1} N_1^\#,\ i \in \Z,
\end{eqnarray}

(3) $S = \{N_1, N_2\}$:
\begin{eqnarray}\label{eq:st_31}
\Lambda(3i) = \p^i N_1,\ 
\Lambda(3i+1) = \p^i N_2,\
\Lambda(3i+2) = \p^{i+1} N_1^\#,\ i \in \Z,
\end{eqnarray}

(4) $S = \{N_0, N_2\}$:
\begin{eqnarray}\label{eq:st_20}
\Lambda(2i) = \p^i N_0,\ 
\Lambda(2i+1) = \p^i N_2,\ i \in \Z,
\end{eqnarray}

(5) $S = \{N_1\}$:
\begin{eqnarray}\label{eq:st_21}
\Lambda(2i) = \p^i N_1,\ 
\Lambda(2i+1) = \p^{i+1} N_1^\#,\ i \in \Z,
\end{eqnarray}

(6) $S = \{N_0\}$:
\begin{eqnarray}\label{eq:st_10}
\Lambda(i) = \p^i N_0,\ i \in \Z,
\end{eqnarray}

(7) $S = \{N_2\}$:
\begin{eqnarray}\label{eq:st_11}
\Lambda(i) = \p^i N_2,\ i \in \Z.
\end{eqnarray}

\begin{rem}\label{rem:strict_2}
By the argument above,
if $\Lambda$ is a strict self-dual $\ri_F$-lattice sequence in
$V$ of period 2, then
we have $[\Lambda(i):\Lambda(i+1)] = q^4$ for all $i\in \Z$.
\end{rem}

\begin{thm}\label{thm:strict}
Let $\pi$ be an irreducible smooth 
representation of $G$ of positive 
level.
Then $\pi$ contains a fundamental skew stratum
$[\Lambda, n, n-1, \beta]$ which satisfies one of the 
following conditions:

(i) $\Lambda$ is a standard strict self-dual $\ri_F$-lattice sequence
and
$(e(\Lambda), n) =1$,

(ii) $\Lambda$ is the strict $\ri_F$-lattice sequence in 
(\ref{eq:st_21})
and $n$ is even.
\end{thm}
\begin{proof}
Recall from \cite{Kariyama}  \S 1.4 that
a $C$-sequence in $V$ is 
a self-dual $\ri_F$-lattice sequence $\Lambda$ in $V$
which satisfies
\begin{enumerate}
\item[$C$(i)] $\Lambda(2i+1) \supsetneq \Lambda(2i+2)$, $i \in \Z$,

\item[$C$(ii)] $e(\Lambda)$ is even and $d(\Lambda)$ is odd.
\end{enumerate}
It follows from  \cite{Kariyama} Proposition 3.1.1 that 
an irreducible smooth representation 
$\pi$ of $G$ of positive level
contains a fundamental skew stratum $[\Lambda, n, n-1, \beta]$ such that
$\Lambda$ is a $C$-sequence and 
$(e(\Lambda), n) = 2$.
After conjugation by some element in $G$,
we may assume $\Lambda$ is a standard $C$-sequence.

We claim that it suffices to find an $\ri_F$-lattice sequence
$\Lambda'$ and an integer $n'$ verifying one of the conditions
in the theorem,
$n/e(\Lambda) = n'/e(\Lambda')$
and
$\rad_{n+1}(\Lambda) \supset \rad_{n'+1}(\Lambda')$.
If this is the case,
then $\pi$ contains a $P_{\Lambda', n'+1}$-fixed vector
and some skew stratum $[\Lambda', n', n'-1, \beta']$.
Proposition~\ref{prop:MP} (iii) says that it is fundamental.

For $\mathcal{L}$ a set of $\ri_F$-lattices,
we write $\mathcal{L}^\# = \{ L^\# \ |\ L \in \mathcal{L}\}$.
Note that 
$\Lambda(2\Z)$ and $\Lambda(2\Z+1)$ are closed 
under the multiplication by elements in $F^\times$,
and 
$\Lambda(2\Z) = \Lambda(2\Z+1)^\#$.

(a) Suppose that 
$\Lambda(2\Z) = \Lambda(2\Z+1)$.
Then we have
$\Lambda(2i) = \Lambda(2i+1)$, $i \in \Z$.
Define an $\ri_F$-lattice sequence 
$\Lambda'$ by
$\Lambda'(i) = \Lambda(2i)$, $i \in \Z$.
Then
$\Lambda'$ is a strict standard self-dual $\ri_F$-lattice
sequence
such that
$e(\Lambda') = e(\Lambda)/2$
and 
$\rad_k(\Lambda') = \rad_{2k-1}(\Lambda) = \rad_{2k}(\Lambda)$, $k \in \Z$.
Putting $n' = n/2$,
we get 
$\rad_{n'+1}(\Lambda') = \rad_{n+1}(\Lambda)$.

(b) 
Suppose that
$\Lambda(2\Z) \cap  \Lambda(2\Z+1) = \emptyset$.
Then 
$\Lambda$ is strict and
$\Lambda(2\Z) \cap  \Lambda(2\Z)^\# = \emptyset$.
Since $N_0^\# = N_0$ and $N_2^\# = \p^{-1}N_2$,
the set
$\Lambda(2\Z)$ is either $\{ \p^m N_1\ |\ m \in \Z\}$
or $\{ \p^m N_1^\#\ |\ m \in \Z\}$,
and hence
$\Lambda(\Z) = \{ \p^m N_1,\ \p^m N_1^\#\ |\ m \in \Z\}$.
Replacing $\Lambda$ with some translate,
we may assume
$\Lambda$ is the $\ri_F$-lattice sequence 
in (\ref{eq:st_21}),
and $[\Lambda, n, n-1,\beta]$ satisfies the condition 
(ii).

(c) 
Suppose that
$\Lambda(2\Z) \neq \Lambda(2\Z+1)$ and
$\Lambda(2\Z) \cap \Lambda(2\Z+1) \neq \emptyset$.
Then we see that 
$\Lambda(\Z)$ contains both $N_1$ and $\p N_1^\#$,
and $\Lambda(2\Z)$ contains
$N_0$ or $N_2$.

Suppose that $\Lambda(2\Z)$ (and hence 
$\Lambda(2\Z+1)$) contains $N_0$.
Replacing $\Lambda$ by some translate,
we may assume 
$\Lambda(0) = \Lambda(1) = N_0$.
Since $N_1$ and $\p N_1^\#$ belong to $\Lambda(\Z)$,
we have $\Lambda(2) = N_1$
and therefore $\Lambda(3)$ 
should be either $N_2$ or $\p N_1^\#$.
If $\Lambda(3) = N_2$,
then $\Lambda(3) = \Lambda(4) = N_2$
since
$N_2 \in \Lambda(2\Z+1)^\# = \Lambda(2\Z)$.
This contradicts to the condition $C$(i).
We therefore have
\begin{eqnarray}\label{eq:b2}
\Lambda(0) = \Lambda(1) = N_0,\ 
\Lambda(2) = N_1,\ \Lambda(3) = \p N_1^\#
\end{eqnarray}
and $e(\Lambda) = 4$.

Similarly,
if $\Lambda(2\Z)$ contains $N_2$,
then
we may assume that
$\Lambda$ satisfies $e(\Lambda) = 4$ and
\begin{eqnarray}\label{eq:b3}
\Lambda(1) = N_1,\ 
\Lambda(2) = \Lambda(3) = N_2,\ 
\Lambda(4) = \p N_1^\#.
\end{eqnarray}

In both cases, we have 
$n = 4m+2$, for $m \geq 0$.
Let $\Lambda'$ be the sequence in (\ref{eq:st_20})
and let
$n' = 2m+1$.
Then it is easy to check that
$\rad_{n+1}(\Lambda) \supset \rad_{n'+1}(\Lambda')$.
This complete the proof.
\end{proof}
\begin{rem}\label{rem:strata}
Replacing $\Lambda$ by some translate,
an irreducible smooth 
representation $\pi$ of $G$ of positive level
contains a fundamental skew stratum 
$[\Lambda, n, n-1, \beta]$
such that
$\Lambda$ and $n$ satisfy one of the conditions listed below:
\begin{eqnarray}
\begin{array}{|c|c|c|c|c|}\hline
\Lambda & e(\Lambda) & d(\Lambda) & (n, e(\Lambda)) & 
e(\Lambda)/(n, e(\Lambda))  \\
\hline
(\ref{eq:st_4})& 4 & 0 & 1 & 4 \\

\hline
(\ref{eq:st_30})& 3 & 0 & 1 & 3  \\
\hline
(\ref{eq:st_31})& 3 & -1 & 1 & 3 \\
\hline
(\ref{eq:st_20})& 2 & 0 & 1 &  2 \\
\hline
(\ref{eq:st_21})& 2 & -1 & 1 & 2 \\
\hline
(\ref{eq:st_10})& 1 & 0 & 1 & 1 \\
\hline
(\ref{eq:st_11})& 1 & -1 & 1 & 1 \\
\hline
(\ref{eq:st_21})& 2 & -1 & 2 & 1 \\
\hline
\end{array}
\end{eqnarray}
\end{rem}
\subsection{Characteristic polynomials}
\label{subsec:characteristic}
Let $\Lambda$ be a strict $\ri_F$-lattice sequence in $V$
and $n$ an integer coprime to $e(\Lambda)$.
Let $[\Lambda, n, n-1, \beta]$ be a stratum in $A$.
We write $e = e(\Lambda)$.
Note that
$\phi_\beta(X)$ 
coincides with the characteristic polynomial of
$y_\beta +\rad_1(\Lambda)$ in 
the $k_F$-algebra $\rad_0(\Lambda)/\rad_1(\Lamdba)$.

As in \cite{HM2} \S 2,
we can form a $\Z/e\Z$-graded $k_F$-vector space
$\overline{\Lambda} = \bigoplus_{i \in \Z/e\Z} \overline{\Lambda}(i)$,
where $\overline{\Lambda}(i) = \Lambda(i)/\Lambda(i+1)$.
There is a natural isomorphism
\[
\rad_k(\Lambda)/\rad_{k+1}(\Lamdba) \simeq 
\bigoplus_{i \in \Z/e\Z}\mathrm{Hom}_{k_F}(\overline{\Lambda}(i), 
\overline{\Lambda}(i+k));
x \mapsto (x_i)_{i \in \Z/e\Z},
\]
given by
\[
x_i:
 v +\Lambda(i+1) \mapsto xv +\Lambda(i+k+1),\
x \in \rad_k(\Lambda),\ v \in \Lambda(i).
\]
For $i \in \Z/e\Z$,
let $\beta_i$ and $y_i$ denote 
the image of $\beta$ and $y_\beta$
in $\mathrm{Hom}_{k_F}(\overline{\Lambda}(i), 
\overline{\Lambda}(i-n))$
and $\mathrm{End}_{k_F}(\overline{\Lambda}(i))$,
respectively.
If we write $\phi_i(X)$ for the characteristic polynomial
of $y_i$,
then we have $\phi_\beta(X) = \prod_{i \in \Z/e\Z} \phi_i(X)$.
Choose $j \in \Z/e\Z$
so that $n_j = \dim_{k_F}(\overline{\Lambda}(j))$ is minimal.
Since $y_i
= \beta_{i-(e-1)n} \cdots \beta_{i-n}\beta_i$,
we have
$\phi_\beta(X) 
= \phi_j(X)^e X^m$,
where $m = \dim_F(V)-en_j$.

Suppose that $[\Lambda, n, n-1, \beta]$ is skew.
Then $\phi_\beta((-1)^eX) = \overline{\phi_\beta(X)}$
since $\sigma(y_\beta) = (-1)^e y_\beta$.
We therefore have
\[
\phi_\beta(X) = 
\left\{
\begin{array}{lll}
(X-a)^4,\ & a\in k_0, & \mathrm{if}\ e = 4,\\
(X-a\e)^3 X,\ & a\in k_0, & \mathrm{if}\ e = 3,\\
(X^2+aX +b)^2,\ & a, b\in k_0, & \mathrm{if}\ e = 2,\\
X^4 +a\e X^3 +bX^2 +c\e X +d,\ & a, b, c, d\in k_0, &
\mathrm{if}\ e = 1.
\end{array}
\right.
\]

Every quadratic polynomial in $k_0[X]$ is reducible in $k_F[X]$.
So, if $e = 2$, then
there are following three cases

(a) $\phi_\beta(X) = (X-\lambda)^2
(X-\overline{\lambda})^2$, $\lambda \in k_F$,
$\lambda \neq \overline{\lambda}$;

(b) $\phi_\beta(X) = (X-u)^2(X-v)^2$, $u, v \in k_0$,
$u \neq v$;

(c) $\phi_\beta(X) = (X-u)^4$, $u \in k_0$.

\section{Representations associated to maximal tori}\label{sec:r_4}
\setcounter{equation}{0}\label{supercuspidal}
The remaining part of this paper is devoted to 
a classification of the 
irreducible smooth representations of $G$
of non-integral level.
In this section,
we classify irreducible representations of $G$ which contain
a semisimple skew stratum 
$[\Lambda, n, n-1, \beta]$ such that
the $G$-centralizer $G_\beta$ of $\beta$ is compact.
We also discuss the \lq\lq intertwining implies conjugacy" property
of such strata.
\subsection{Construction of supercuspidal representations}
We first recall the construction of supercuspidal representations
from \cite{St2} and \cite{St1}.

Let $[\Lambda, n, n-1, \beta]$ be a semisimple skew stratum in $A$.
We write $G_\beta$ for the $G$-centralizer of $\beta$.
By \cite{St2} Theorem 4.6, we get
\begin{eqnarray}
I_G[\Lambda, n, [n/2], \beta]
= P_{\Lambda, [(n+1)/2]} G_\beta P_{\Lambda, [(n+1)/2]}.
\end{eqnarray}

Suppose that $G_\beta$ is compact.
Then $G_\beta$ lies in $P_{\Lambda, 0}$,
and hence
$I_G[\Lambda, n, [n/2], \beta]
= G_\beta P_{\Lambda, [(n+1)/2]}$.
This implies that 
the intertwining of $(P_{\Lambda, [n/2]+1}, \psi_\beta)$
coincides with the normalizer of $(P_{\Lambda, [n/2]+1}, \psi_\beta)$.
So we obtain the following:
\begin{prop}\label{prop:classification}
Let $[\Lambda, n, n-1, \beta]$ be a semisimple skew stratum 
such that $G_\beta$ is compact.
Set $J = G_\beta P_{\Lambda, [(n+1)/2]}$.
Then 
the map $\rho \mapsto \mathrm{Ind}_J^G \rho$ 
gives a bijection between 
$\Irr(J)^{(P_{\Lambda, [n/2]+1}, \psi_\beta)}$ and
$\Irr(G)^{(P_{\Lambda, [n/2]+1}, \psi_\beta)}$.
\end{prop}
\begin{proof}
For $\rho \in \Irr(J)^{(P_{\Lambda, [n/2]+1}, \psi_\beta)}$,
the restriction of $\rho$ to $P_{\Lambda, [n/2]+1}$
is a multiple of $\psi_\beta$ since
$J$ normalizes $\psi_\beta$.
If $g \in G$ intertwines $\rho$,
then $g$ intertwines $\psi_\beta$,
and hence $g$ lies in $J$.
Thus, $\mathrm{Ind}_J^G \rho$ is irreducible and supercuspidal.

The surjectivity of this map follows from Frobenius reciprocity.
Let $\rho_i \in \Irr(J)^{(P_{\Lambda, [n/2]+1}, \psi_\beta)}$,
$i = 1, 2$.
Suppose that 
$\mathrm{Ind}_J^G \rho_1$ is isomorphic to 
$\mathrm{Ind}_J^G \rho_2$.
Then there is $g \in G$ which intertwines $\rho_1$ and $\rho_2$.
Restricting these to $P_{\Lambda, [n/2]+1}$,
we get $g \in J$, and $\rho_1 \simeq \rho_2^g \simeq \rho_2$.
This implies the injectivity.
\end{proof}

We shall give a description of $\Irr(J)^{(P_{\Lambda, [n/2]+1}, \psi_\beta)}$ when 
$G_\beta$ is a maximal torus in $G$.
In this case, the group $G_\beta$ is abelian.
Put $H^1 = (G_\beta \cap P_{\Lambda, 1})P_{\Lambda, [n/2]+1}$
and $J^1 = (G_\beta \cap P_{\Lambda, 1})P_{\Lambda, [(n+1)/2]}$.
Let $\sigma \in \Irr(J)^{(P_{\Lambda, [n/2]+1}, \psi_\beta)}$.
Since $H^1/P_{\Lambda, [n/2]+1} \simeq (G_\beta \cap 
P_{\Lambda, 1})/(G_\beta \cap P_{\Lambda, [n/2]+1})$,
$\sigma$ contains some extension $\theta$ of $\psi_\beta$ 
to $H^1$.
Applying the proof of \cite{St1} Proposition 4.1,
we see that there exists a unique irreducible representation 
$\eta_\theta$ of $J^1$ containing $\theta$.
Moreover,
the restriction of  $\eta_\theta$ to $H^1$ is the 
$[J^1: H^1]^{1/2}$-multiple of $\theta$.
The group
$J/J^1 \simeq G_\beta/G_{\beta}\cap P_{\Lambda,1}$
is cyclic with order coprime to $[J^1: H^1]^{1/2}$
since $[J^1: H^1]^{1/2}$ is a power of $q$.
Hence $\eta_\theta$ can extend to an irreducible representation
$\kappa_\theta$ of  $J$.
It is easy to observe that
$\Irr(J)^{(J^1, \eta_\theta)} = 
\{ \chi \otimes \kappa_\theta\ |\ \chi \in 
(G_\beta/G_{\beta}\cap P_{\Lambda,1})^\wedge \}$,
and $\chi \otimes \kappa_\theta \simeq \chi' \otimes 
\kappa_\theta$
if and only if $\chi = \chi'$,
for $\chi, \chi' \in (G_\beta/G_{\beta}\cap P_{\Lambda,1})^\wedge$.
In particular,
the restriction of every element in $\Irr(J)^{(J^1, \eta_\theta)}$
to $H^1$ is a multiple of $\theta$.
We conclude the following:
\begin{prop}\label{prop:description}
With the notation as above,
the set $\Irr(J)^{(P_{\Lambda, [n/2]+1}, \psi_\beta)}$
is a disjoint union $\bigcup_{\theta}\Irr(J)^{(H^1, \theta)}$,
where $\theta$ runs over the extensions of $\psi_\beta$ to $H^1$.
We have
$\Irr(J)^{(H^1, \theta)} = 
\{ \chi \otimes \kappa_\theta\ |\ \chi \in 
(G_\beta/G_{\beta}\cap P_{\Lambda,1})^\wedge \}$.
\end{prop}

\subsection{Representations of level $n/4$}\label{subsec:max}
In this section, we fix a positive integer $n$ coprime to 4
and 
the $\ri_F$-lattice sequence $\Lambda$
in
(\ref{eq:st_4}).
By Remark~\ref{rem:strata},
an irreducible smooth representation of $G$ of level $n/4$
contains a fundamental skew stratum $[\Lambda, n, [n/2], \beta]$.

Let $\Lambda'$ be a
strict self-dual $\ri_F$-lattice sequence in $V$
of period 4
and $[\Lambda', n, n-1, \beta]$ a fundamental skew stratum.
Proposition~\ref{prop:pure} says that
this is a skew simple stratum and 
$E_\beta$ is a totally ramified extension of degree 4 over $F$,
where $E_\beta = F[\beta]$.
Propositions~\ref{prop:classification} and 
\ref{prop:description} give a classification of 
the irreducible smooth representations of $G$ which contain
$[\Lambda', n, [n/2], \beta]$.

We give a necessary and sufficient condition
when two  such strata occur in a common irreducible
representation of $G$.
\begin{thm}\label{thm:conj_4}
Let  
$\Lambda$ be the $\ri_F$-lattice sequence 
in (\ref{eq:st_4}) and let
$n$ be a positive integer coprime to $4$.
Let 
$[\Lambda, n, [n/2], \beta]$ and
$[\Lambda, n, [n/2], \gamma]$ be fundamental skew strata.
Suppose that there is an irreducible smooth representation of $G$
containing $[\Lambda, n, [n/2], \beta]$ and
$[\Lambda, n, [n/2], \gamma]$.
Then 
$(P_{\Lambda, [n/2]+1}, \psi_\gamma)$ is 
a $P_{\Lambda, 0}$-conjugate of $(P_{\Lambda, [n/2]+1}, \psi_\beta)$.
\end{thm}
\begin{proof}
By Proposition~\ref{prop:MP} (i), there is $g \in G$ such that
$(\beta + \g_{\Lambda, -[n/2]})\cap 
\Ad(g)(\gamma + \g_{\Lambda, -[n/2]})\neq \emptyset$.
Take $\delta \in (\beta + \g_{\Lambda, -[n/2]})\cap 
\Ad(g)(\gamma + \g_{\Lambda, -[n/2]})$.
We obtain fundamental skew strata 
$[\Lambda, n, n-1, \delta]$ and
$[g\Lambda, n, n-1, \delta]$.
We see that $\Lambda$ and $g\Lambda$ are strict $\ri_{E_\delta}$-lattice sequence in the one dimension $E_\delta$-space $V$,
and hence $\Lambda$ is a translate of $g\Lambda$.
Since
$d(g\Lambda) = d(\Lambda)$,
we get $\Lambda = g\Lambda$ and $g \in P_{\Lambda, 0}$.
This completes the proof.
\end{proof}

\subsection{Representations of level $n/3$}\label{subsec:r_3}
\label{sec:r_3}
We fix a positive integer $n$ coprime to 3.
By Remark~\ref{rem:strata},
an irreducible smooth representation of $G$ of level $n/3$
contains a fundamental skew stratum $[\Lambda, n, [n/2], \beta]$
such that
$\Lambda$ is the $\ri_F$-lattice sequence in
(\ref{eq:st_30}) or (\ref{eq:st_31}).

Let 
$\Lambda$ be any strict self-dual $\ri_F$-lattice sequence in $V$
with $e(\Lambda) = 3$
and 
let 
$[\Lambda, n, [n/2], \beta]$ be a fundamental skew stratum.
As seen in \S \ref{subsec:characteristic},
we have 
$\phi_\beta(X) = (X-a\e)^3 X$,
for $a \in k_0^\times$.
Using Hensel's Lemma,
we can lift this to
$\Phi_\beta(X) = f_a(X) f_0(X)$
where $f_a(X)$, $f_0(X)$ are monic,
$f_a(X) \bmod{\mi_F} = (X-a\e)^3$
and $f_0(X)  \bmod{\mi_F} =  X$.
Applying the argument in the proof of \cite{St3} Theorem 4.4,
if we put
$V^{a} = \ker f_a(y_\beta)$ and  $V^{0} = \ker f_0(y_\beta)$,
then $[\Lambda, n, [n/2], \beta]$ is split 
with respect to 
 $V = V^{a} \bot V^{0}$.
As usual, for $b \in \{a, 0\}$,
we write $\beta_b = \beta|_{V^b}$ and 
$\Lambda^b(i) = \Lambda(i) \cap V^b$, $i \in \Z$.

\begin{lem}\label{lem:over3}
The form $(V^0, f|_{V^0})$ represents 1
if and only if $d(\Lambda)$ is even.
\end{lem}
\begin{proof}
Suppose that $d(\Lambda) = 2k$, $k \in \Z$.
Then the dual lattice of $\Lambda^{0}(k)$ 
in $(V^{0}, f|_{V^{0}})$ is itself.
Since $\dim_F V^{0} = 1$, we see that
the form $f|_{V^{0}}$ represents 1.

If $d(\Lambda) = 2k+1$, $k \in \Z$,
then the dual of $\Lambda^{0}(k + 2)$ in $(V^{0}, f|_{V^{0}})$ is
$\p^{-1}\Lambda^{0}(k +2)$.
The form $f|_{V^{0}}$ represents $\p$,
and does not represent 1 since 
$F$ is unramified over $F_0$.
\end{proof}

\begin{prop}\label{prop:disjoint_3}
Let 
$\Lambda$ and $\Lambda'$ be the $\ri_F$-lattice sequence defined by 
(\ref{eq:st_30}) and (\ref{eq:st_31}), respectively.
Let 
$[\Lambda, n, [n/2], \beta]$ and $[\Lambda', n, [n/2], \gamma]$
be fundamental skew strata.
Then, there are no irreducible smooth representations of $G$
which contain $[\Lambda, n, [n/2], \beta]$ and 
$[\Lambda, n,  [n/2], \gamma]$.
\end{prop}
\begin{proof}
Suppose that 
there exists an irreducible smooth representation of $G$
containing $[\Lambda, n, [n/2], \beta]$ and 
$[\Lambda, n,  [n/2], \gamma]$.
By Proposition~\ref{thm:MP} (i), 
we can take 
$\delta \in (\beta + \g_{\Lambda, -[n/2]})\cap 
(\Ad(g)\gamma  + \g_{g\Lambda', -[n/2]})$,
for some $g \in G$.

Applying the above construction of a splitting to 
fundamental skew strata
$[\Lambda, n, n-1, \delta]$ and 
$[g \Lambda', n, n-1, \delta]$,
it follows from Lemma~\ref{lem:over3} that
$d(\Lambda) \equiv d(g\Lambda') \pmod{2}$, which 
contradicts the choice of $\Lambda$ and $\Lambda'$.
This completes the proof.
\end{proof}

Due to \cite{K4} Lemma 3.6, 
$\Lambda^{a}$ is a strict $\ri_F$-lattice sequence 
of period 3 in 
the 3-dimensional $F$-space $V^{a}$,
and 
by Proposition~\ref{prop:pure},
$[\Lambda^{a}, n, n-1, \beta_a]$ is simple.
Let $E_a = F[\beta_a]$.
Since $\Lambda^a$ is a strict $\ri_{E_a}$-lattice sequence 
in the one-dimensional $E_a$-space $V^a$,
the integer $d(\Lambda)$ determines
$\Lambda^a$ uniquely.
Similarly,
$\Lambda^0$ is an $\ri_F$-lattice sequence in 
the one-dimensional $F$-space $V^0$ of period 3,
so it is also determined by $d(\Lambda)$.
We therefore see that given such a stratum 
$[\Lambda, n, n-1, \beta]$,
we can reconstruct $\Lambda$ from $d(\Lambda)$.

We get the following theorem,
imitating the proof of Theorem~\ref{thm:conj_4}.
\begin{thm}
Let 
$\Lambda$ be the $\ri_F$-lattice sequence in
(\ref{eq:st_30}) or (\ref{eq:st_31}),
 and let
$n$ be a positive integer coprime to $3$.
Let 
$[\Lambda, n, [n/2], \beta]$ and
$[\Lambda, n, [n/2], \gamma]$ be fundamental skew strata.
Suppose that there exists an irreducible smooth representation
of $G$ which contains $[\Lambda, n, [n/2], \beta]$
and $[\Lambda, n, [n/2], \gamma]$.
Then
$(P_{\Lambda, [n/2]+1}, \psi_\gamma)$ is 
a $P_{\Lambda, 0}$-conjugate of $(P_{\Lambda, [n/2]+1}, \psi_\beta)$.
\end{thm}

Let $[\Lambda, n, [n/2], \beta]$ be a skew stratum
as above.
The restriction $[\Lambda, n, n-1, \beta]$ is not 
semisimple only if $\beta_2 \neq 0$
since $\phi_{\beta_2}(X) = X$.
Even in this case,
$\beta_2$ is scalar and does not matter 
in computing the intertwining.
Hence 
we can classify 
the irreducible smooth representations of $G$ which contain
$[\Lambda, n, [n/2], \beta]$,
using
Propositions~\ref{prop:classification} and 
\ref{prop:description}.

\subsection{Representations of half-integral level}
\label{subsec:half}
We fix a positive odd integer $n$.
Let $\pi$ be an irreducible smooth representation of $G$
of level $n/2$ and of characteristic polynomial
\begin{eqnarray}\label{eq:char_2}
(X-a)^2 (X-b)^2,\ 
\mathrm{for}\ a\in k_0^\times,\
b \in k_0,\ a\neq b.
\end{eqnarray}
Then by Remark~\ref{rem:strata},
$\pi$ contains a fundamental skew stratum
$[\Lambda, n, [n/2], \beta]$ such that
$\Lambda$ is the sequence in (\ref{eq:st_20}) or (\ref{eq:st_21})
and $\phi_\beta(X)$ has the form (\ref{eq:char_2}).

Let $\Lambda$ be any strict self-dual $\ri_F$-lattice sequence
in $V$ with $e(\Lambda) = 2$
and let
$[\Lambda, n, n-1, \beta]$ be a fundamental skew stratum
such that $\phi_\beta(X)$ has the form (\ref{eq:char_2}).
As in \S \ref{sec:r_3},
by Hensel's Lemma,
we can lift  $\phi_\beta(X)$ to
$\Phi_\beta(X) = f_{a}(X) f_{b}(X)$
where, for $c \in \{ a, b \}$,
$f_{c}(X)$ is monic and
$f_{c}(X) \bmod{\mi_F} = (X-c)^2$.
As in the proof of \cite{St3} Theorem 4.4,
if we put
$V^{a} = \ker f_{a}(y_\beta)$ and
$V^{b} = \ker f_{b}(y_\beta)$,
then $[\Lambda, n, [n/2], \beta]$ is split with respect to 
$V = V^{a} \bot V^{b}$.

For $c \in \{ a, b\}$,
we write $\beta_c = \beta|_{V^{c}}$
and $\Lambda^c(i) = \Lambda(i) \cap V^c$, $i \in \Z$.
It follows from \cite{K4} Lemma 3.6 that
$\Lambda^{a}$
is a strict $\ri_F$-lattice sequence in $V^{a}$ of
period 2,
and hence $[\Lambda^{a}(i):\Lambda^{a}(i+1)]
= q^2$ for all $i \in \Z$.
By Remark~\ref{rem:strict_2},
$\Lambda^{b}$ is also strict.
\begin{lem}\label{lem:st_h_int}
(i) 
If $d(\Lambda)$ is even,
then 
$(V^{a}, f|_{V^{a}})$ and $(V^{b},
f|_{V^{b}})$ are 
isotropic.

(ii)
If $d(\Lambda)$ is odd,
then 
$(V^{a}, f|_{V^{a}})$ and $(V^{b},
f|_{V^{b}})$ are 
anisotropic.
\end{lem}
\begin{proof}
The assertion follows from
Lemma~\ref{lem:anisotropic}.
\end{proof}

\begin{prop}
Let 
$\Lambda$ and $\Lambda'$ be the $\ri_F$-lattice sequences in
(\ref{eq:st_20}) and (\ref{eq:st_21}), respectively.
Let $[\Lambda, n, [n/2], \beta]$ and $[\Lambda', n, [n/2], \gamma]$ be fundamental skew strata.
Suppose 
$\phi_\beta(X) = \phi_\gamma(X) = (X-a)^2 (X-b)^2$,
for $a, b \in k_0$ such that $a \neq b$.
Then there are no irreducible smooth representations of $G$
which contain $[\Lambda, n, [n/2], \beta]$ 
and $[\Lambda, n, [n/2], \gamma]$.
\end{prop}
\begin{proof}
The assertion follows 
from Lemma~\ref{lem:st_h_int}
as in the proof of  Proposition~\ref{prop:disjoint_3}.
\end{proof}

For $c \in \{ a, b\}$, 
we put $E_c = F[\beta_c]$.
Proposition~\ref{prop:pure} says 
that
$[\Lambda^{a}, n, n-1, \beta_a]$ is simple
and $\Lambda^a$ is a strict $\ri_{E_a}$-lattice sequence 
in the one-dimensional $E_a$-space $V^a$.
An analog holds for $[\Lambda^{b}, n, n-1, \beta_b]$
if $b \neq 0$.
In this case, $d(\Lambda)$ determines $\Lambda$.

Suppose that $d(\Lambda)$ is odd.
Lemma~\ref{lem:st_h_int} implies the uniqueness
of $\Lambda^c$ in the anisotropic space $(V^c, f|_{V^c})$
up to translation.
Hence $d(\Lambda)$ determines $\Lambda$.
\begin{thm}\label{thm:sch}
Let 
$\Lambda$ be the $\ri_F$-lattice sequence in
(\ref{eq:st_20}) or (\ref{eq:st_21}),
and let
$n$ a positive odd integer.
Let 
$[\Lambda, n, [n/2], \beta]$ and
$[\Lambda, n, [n/2], \gamma]$ be skew strata
such that
$\phi_\beta(X) = \phi_\gamma(X) = 
(X-a)^2 (X-b)^2$, $a \neq b$.
Assume that $ab \neq 0$ or
$\Lambda$ is defined by (\ref{eq:st_21}).
Suppose that there exists an irreducible smooth representation of $G$
containing $[\Lambda, n, [n/2], \beta]$
and $[\Lambda, n, [n/2], \gamma]$.
Then
$(P_{\Lambda, [n/2]+1}, \psi_\gamma)$ is 
a $P_{\Lambda, 0}$-conjugate of $(P_{\Lambda, [n/2]+1}, \psi_\beta)$.
\end{thm}
\begin{proof}
This is exactly as in the proof of Theorem~\ref{thm:conj_4}.
\end{proof}

We shall describe the set $\Irr(G)^{(P_{\Lambda, [n/2]+1}, \psi_\beta)}$ under the assumption of Theorem~\ref{thm:sch}.
If  $ab \neq 0$,
then we can apply
Propositions~\ref{prop:classification} and 
\ref{prop:description}.

Suppose that $\Lambda$ is defined by 
(\ref{eq:st_21}) and $b = 0$.
We write 
$G_c= G\cap \mathrm{End}_F(V^{c})$,
for $c \in \{ a, 0\}$.
Let $G_{\beta_a}$ denote the $G_a$-centralizer of $\beta_a$
and 
$I_0$ the formal intertwining of $[\Lambda^{0},
n, [n/2], \beta_0]$ in $G_0$.
Since $G_{0}$ coincides with
$P_{\Lambda^0, 0}$,
$I_0$ is just 
the normalizer of $(P_{\Lambda, [n/2]+1}\cap G_0, \psi_\beta
|_{P_{\Lambda, [n/2]+1}\cap G_0})$
in $G_0$.
By
\cite{St2} Theorem 4.6 and Corollary 4.14, we have
\begin{eqnarray}
I_G[\Lambda, n, [n/2], \beta]
= G' P_{\Lambda, [n/2]+1},
\end{eqnarray}
where $G' = G_{\beta_a} \times I_0$.
Therefore,
we get a natural isomorphism of Hecke algebras
\begin{eqnarray}
\He(G//P_{\Lambda, [n/2]+1}, \psi_\beta)
\simeq
\He(G'//P', \psi'),
\end{eqnarray}
where $P' = G' \cap P_{\Lambda, [n/2]+1}$ and
$\psi' = \psi_\beta|_{P'}$.
Via this isomorphism,
we can identify
$\mathrm{Irr}(G)^{(P_{\Lambda, [n/2]+1}, \psi_\beta)}$ with
$\mathrm{Irr}(G')^{(P', \psi')}$.
Moreover, 
every irreducible representation of $G$ containing $[\Lambda, n, [n/2], \beta]$
is supercuspidal since $G'$ is compact.

\section{Hecke algebra isomorphisms}\label{sec:h_int}
\setcounter{equation}{0}\label{Hecke}
In this section, 
we complete our classification of the
irreducible smooth representations of $G$ of non-integral level.
The remaining representations we need to consider 
are those of half-integral level.
We fix a positive integer $m$ and put $n = 2m-1$.
By Proposition~\ref{thm:MP},
the set of equivalence classes of irreducible smooth representations
of $G$ of level $n/2$ is roughly
classified according to those characteristic polynomials.
Let $\Lambda$ be a strict self-dual $\ri_F$-lattice sequence
of period 2.
As claimed in \S~\ref{subsec:characteristic},
the characteristic polynomial of a fundamental skew stratum
$[\Lambda, n, n-1, \beta]$ has one of the following forms:

(\ref{sec:h_int}a) $(X -\lambda)^2(X-\overline{\lambda})^2$, $\lambda\in k_F$,
$\lambda \neq \overline{\lambda}$;

(\ref{sec:h_int}b)
$(X -u)^2(X-v)^2$, $u \in k_0^\times$,
$v \in k_0$,
$u \neq v$;

(\ref{sec:h_int}c)
$(X -u)^4$, $u \in k_0^\times$.

\noindent
We have classified the representations of case (\ref{sec:h_int}b)
in \S \ref{subsec:half},
except when

(\ref{sec:h_int}d)
$(X-u)^2 X^2$, $u \in k_0^\times$ and 
$\Lambda$ is defined by (\ref{eq:st_20}).
\subsection{Normalization of $\beta$}\label{subsec:norm}
Let $\Lambda$ and $\Lambda'$
be the $\ri_F$-lattice sequences in (\ref{eq:st_20})
and (\ref{eq:st_21}), respectively.
By Remark~\ref{rem:strata},
an irreducible smooth representation of $G$ of level $n/2$
contains a fundamental skew stratum
of the form $[\Lambda, n, n-1,\beta]$
or $[\Lambda', n, n-1, \beta']$.

The filtrations
$\{\seqrad_k(\Lambda)\}_{k \in \Z}$  and 
$\{\rad_k(\Lambda')\}_{k \in \Z}$
are given by the following:
\begin{eqnarray}\label{eq:fil_20}
\seqrad_0(\Lambda) = 
\left(
\begin{array}{cc|cc}
\ri_F & \ri_F & \ri_F & \ri_F \\
\ri_F & \ri_F & \ri_F & \ri_F \\ \hline
\mi_F & \mi_F & \ri_F & \ri_F \\
\mi_F & \mi_F & \ri_F & \ri_F \\
\end{array}
\right),\
\seqrad_1(\Lambda) = 
\left(
\begin{array}{cc|cc}
\mi_F & \mi_F & \ri_F & \ri_F \\
\mi_F & \mi_F & \ri_F & \ri_F \\ \hline
\mi_F & \mi_F & \mi_F & \mi_F \\
\mi_F & \mi_F & \mi_F & \mi_F \\
\end{array}
\right),
\end{eqnarray}
\begin{eqnarray}
\rad_0(\Lambda')
= 
\left(
\begin{array}{c|cc|c}
\ri_F & \ri_F & \ri_F & \mi_F^{-1}\\ \hline
\mi_F & \ri_F & \ri_F & \ri_F \\
\mi_F & \ri_F & \ri_F & \ri_F \\ \hline
\mi_F & \mi_F & \mi_F & \ri_F 
\end{array}
\right),\
\rad_1(\Lambda')
= 
\left(
\begin{array}{c|cc|c}
\mi_F & \ri_F & \ri_F & \ri_F\\ \hline
\mi_F & \mi_F & \mi_F & \ri_F \\
\mi_F & \mi_F & \mi_F & \ri_F \\ \hline
\mi_F^2 & \mi_F & \mi_F & \mi_F 
\end{array}
\right).
\end{eqnarray}
Up to equivalence classes of skew strata,
we can take
$\beta \in \g_{\Lambda, -n}$
and $\beta' \in \g_{\Lambda', -n}$
to be
\begin{eqnarray}
& \beta = \varpi^{-m}\left(\begin{array}{cccc}
0 & 0 & C & a\sqrt{\varepsilon} \\
 0 & 0 & b\sqrt{\varepsilon} & -\overline{C} \\ 
\varpi D & \varpi c\sqrt{\varepsilon} & 0 & 0 \\
\varpi d\sqrt{\varepsilon} & -\varpi \overline{D} & 0 & 0 \\
\end{array}
\right),\ C, D \in \ri_F,\ a,b,c,d \in \ri_0;
\end{eqnarray}
\begin{eqnarray}
& \beta' = \p^{-m}
\left(
\begin{array}{cccc}
0 & Z & C & 0\\ 
\p M & 0 & 0 & -\overline{C}\\
\p D & 0 & 0 & -\overline{Z}\\
0 & -\p \overline{D} & -\p \overline{M} & 0
\end{array}
\right),\ C, D, M, Z \in \ri_F.
\label{eq:h2}
\end{eqnarray}

\begin{lem}\label{lem:lambda}
Let $[\Lambda', n, n-1, \beta']$ be a fundamental skew stratum 
whose characteristic polynomial is of the form (\ref{sec:h_int}a)
or (\ref{sec:h_int}c).
Suppose that an irreducible smooth 
representation $\pi$ of $G$ contains 
$[\Lambda', n, n-1, \beta']$.
Then there exists
a fundamental skew stratum $[\Lambda, n, n-1, \beta]$
which occurs in $\pi$.
\end{lem}
\begin{proof}
Take $\beta'$ as in (\ref{eq:h2}).
By assumption,
after conjugation by an element in $P_{\Lambda', 0}$,
we may assume $y_\beta$ is upper triangular.
By replacing $\beta$ with a $P_{\Lambda', 0}$-conjugate again,
we may assume that $Z \in \mi_F$.
Then we have $\beta' + \rad_{1-n}(\Lambda') \subset 
\seqrad_{-n}(\Lambda)$.
Thus the assertion follows 
from the standard argument.
\end{proof}

By Lemma~\ref{lem:lambda},
it suffices to consider irreducible smooth representations of $G$ containing a fundamental skew stratum $[\Lambda,
n, n-1, \beta]$ such that 
$\Lambda$ 
is the $\ri_F$-lattice sequences in (\ref{eq:st_20})
and 
$\phi_\beta(X)$ has the form
(\ref{sec:h_int}a), (\ref{sec:h_int}c) or (\ref{sec:h_int}d).
We shall normalize $\beta$
as in \cite{GSp4} Case (5.1b).
Set
\[
X = 
\left(
\begin{array}{cc}
C & a\e\\
b\e & -\overline{C}
\end{array}\right),\ 
Z 
= \left(
\begin{array}{cc}
D & c\e \\
d\e & -\overline{D}
\end{array}
\right),\
S = \left(
\begin{array}{cc}
0 & 1\\
1 & 0
\end{array}
\right)
\]
and abbreviate
\[
\beta = \p^{-m}
\left(
\begin{array}{cc}
0 & X \\
\p Z & 0
\end{array}
\right).
\]
Then we have
\[
y_\beta = \left(
\begin{array}{cc}
XZ & 0\\
0 & ZX
\end{array}
\right).
\]
For $Y \in GL_2(\ri_F)$, we define an element
$g = g(Y)$ in $P_{\Lambda, 0}$ by 
\[
g = \left(
\begin{array}{cc}
Y & 0\\
0 & S{}^t\overline{Y}^{-1}S
\end{array}
\right).
\]
Without loss,
we can replace
$[\Lambda, n, n-1, \beta]$ by
$[\Lambda, n, n-1, Ad(g)(\beta)]$.
Observe that
\[
Ad(g)(\beta) = 
\p^{-m}\left(
\begin{array}{cc}
0 & YXS{}^t\overline{Y}S\\
\p S{}^t\overline{Y}^{-1}SZY^{-1} & 0
\end{array}
\right),
\]
and
\[
Ad(g)(y_\beta) = \left(
\begin{array}{cc}
Ad(Y)(XZ) & 0\\
0 & Ad(S{}^t\overline{Y}^{-1}S)(ZX)
\end{array}
\right).
\]
We will confuse elements in $k_F$ with those in $\ri_F$.

\noindent
Case (\ref{sec:h_int}a):
By Hensel's lemma,
we may assume that
the characteristic polynomial of $XZ$
has the form $(X-\lambda)(X-\overline{\lambda})$, for 
$\lambda \in \ri_F$ such that
$\lambda \not\equiv \overline{\lambda} \pmod{\mi_F}$.
Replacing $\beta$ by a $P_{\Lambda, 0}$-conjugate,
we may assume 
\[
XZ = 
\left(
\begin{array}{cc}
\lambda & 0\\
0 & \overline{\lambda}
\end{array}
\right). 
\]
The spaces $Fe_{1} \oplus Fe_{3}$ and $Fe_2 \oplus Fe_4$
are the kernels of
$y_\beta -\lambda$ 
and $y_\beta -\overline{\lambda}$, respectively,
so that
these spaces are
stable under the action of $\beta$.
Thus $X$ and $Z$ are diagonal.
After conjugation by a diagonal matrix in $P_{\Lambda, 0}$,
we may assume 
\begin{eqnarray}
X = \left(
\begin{array}{cc}
1 & 0\\
0 & -1
\end{array}
\right),\
Z = \left(
\begin{array}{cc}
\lambda & 0\\
0 & -\overline{\lambda}
\end{array}
\right).
\end{eqnarray}

\noindent
Case (\ref{sec:h_int}c):
After $P_{\Lambda, 0}$-conjugation,
we may assume
$XZ$ is upper triangular modulo $\mi_F$.
After $g(Y)$-conjugation by some upper triangular 
$Y \in GL_N(\ri_F)$,
we may also assume 
$X$ is antidiagonal or upper triangular modulo $\mi_F$.

If $X$ is antidiagonal,
then $Z$ is also antidiagonal hence
$XZ$ is scalar.
We can take $X$ to be upper triangular
by replacing $\beta$ with 
a $P_{\Lambda, 0}$-conjugate.
So we can always assume $X$ is upper triangular.
Then $Z$ must be upper triangular as well.
Replacing $\beta$ by some $P_{\Lambda, 0}$-conjugate,
we can assume
\begin{eqnarray}
X = \left(
\begin{array}{cc}
1 & 0\\
0 & -1
\end{array}
\right),\
Z = \left(
\begin{array}{cc}
u & c\e\\
0 & -u
\end{array}
\right),\
u \in \ri_0^\times, c \in \ri_F.
\end{eqnarray}

\noindent
Case (\ref{sec:h_int}d):
As in the case (\ref{sec:h_int}a),
we may assume 
\[
XZ = 
\left(
\begin{array}{cc}
u & 0\\
0 & v
\end{array}
\right),\ u \in \ri_0^\times, v \in \mi_0,
\]
so that $X$ and $Z$ are both anti-diagonal.
Replacing $\beta$ by a $P_{\Lambda, 0}$-conjugate,
we may assume that
$\p^m \beta$ is one of the following elements:
\begin{eqnarray*} 
\left(
\begin{array}{cccc}
0 & 0 & 0 & \e\\
0 & 0 & 0 & 0\\
0 & 0 & 0 & 0\\
\p d \e & 0 & 0 & 0
\end{array}
\right),\
\left(
\begin{array}{cccc}
0 & 0 & 0 & \e\\
0 & 0 & \e & 0\\
0 & 0 & 0 & 0\\
\p d \e & 0 & 0 & 0
\end{array}
\right),\
\left(
\begin{array}{cccc}
0 & 0 & 0 & \e\\
0 & 0 & 0 & 0\\
0 & \p \e & 0 & 0\\
\p d \e & 0 & 0 & 0
\end{array}
\right),
\end{eqnarray*}
where $d = 
u \varepsilon^{-1} \in \ri_0^\times$.
\begin{rem}
In case (\ref{sec:h_int}a),
we can swap $\lambda$ for $\overline{\lambda}$.
\end{rem}

\subsection{Case (\ref{sec:h_int}c)}\label{u11_E}
Let 
$\Lambda'$ be the $\ri_F$-lattice sequence defined by
(\ref{eq:st_20}).
We fix a positive integer $m$ and $\alpha \in \ri_0^\times$. 
Put $n' = 2m-1$.
For $c \in \ri_0$, set
\begin{eqnarray}
\beta(c) = \p^{-m}
\left(
\begin{array}{cccc}
0 & 0 & 1 & 0\\
0 & 0 & 0 & -1\\ 
\p\alpha & \p c\e & 0 & 0\\
0 & -\p\alpha & 0 & 0
\end{array}
\right) \in \g_{\Lambda', -n'}.
\end{eqnarray}

Define a self-dual $\ri_F$-lattice sequence $\Lambda$ in $V$ with $e(\Lambda) = 8$ and $d(\Lambda) = 2$
by
\begin{eqnarray*}
& \Lambda(0) = \Lambda(1) = \Lambda(2) = N_0,\
\Lambda(3) = N_1,\
\Lambda(4) = \Lambda(5) = \Lambda(6) = N_2,\
\Lambda(7) = \p N_1^\#,
\end{eqnarray*}
\begin{eqnarray}\label{eq:L2}
\Lambda(i + 8k) = \p^k \Lambda(i),\ 
0 \leq i\leq 7,\  k \in \Z.
\end{eqnarray}
Then
the induced filtration
$\{\rad_n(\Lambda)\}_{n \in \Z}$ 
is given by
\[
\rad_0(\Lambda) = \left(
\begin{array}{cc|cc}
\ri_F & \ri_F & \ri_F & \ri_F\\
\mi_F & \ri_F & \ri_F & \ri_F\\ \hline
\mi_F & \mi_F & \ri_F & \ri_F\\
\mi_F & \mi_F & \mi_F & \ri_F\\
\end{array}
\right),\
\rad_1(\Lambda) = \left(
\begin{array}{cc|cc}
\mi_F & \ri_F & \ri_F & \ri_F\\
\mi_F & \mi_F & \ri_F & \ri_F\\ \hline
\mi_F & \mi_F & \mi_F & \ri_F\\
\mi_F & \mi_F & \mi_F & \mi_F\\
\end{array}
\right),
\]
\[
\rad_2(\Lambda)= \rad_3(\Lambda) = \left(
\begin{array}{cc|cc}
\mi_F & \mi_F & \ri_F & \ri_F\\
\mi_F & \mi_F & \ri_F & \ri_F\\ \hline
\mi_F & \mi_F & \mi_F & \mi_F\\
\mi_F & \mi_F & \mi_F & \mi_F\\
\end{array}
\right),\
\rad_4(\Lambda)= \left(
\begin{array}{cc|cc}
\mi_F & \mi_F & \ri_F & \ri_F\\
\mi_F & \mi_F & \mi_F & \ri_F\\ \hline
\mi_F & \mi_F & \mi_F & \mi_F\\
\mi_F^2 & \mi_F & \mi_F & \mi_F\\
\end{array}
\right),
\]
\begin{eqnarray}\label{eq:E_fil}
\rad_5(\Lambda)= \left(
\begin{array}{cc|cc}
\mi_F & \mi_F & \mi_F & \ri_F\\
\mi_F & \mi_F & \mi_F & \mi_F\\ \hline
\mi_F^2 & \mi_F & \mi_F & \mi_F\\
\mi_F^2 & \mi_F^2 & \mi_F & \mi_F\\
\end{array}
\right),\
\rad_6(\Lambda)=  \rad_7(\Lambda) = \left(
\begin{array}{cc|cc}
\mi_F & \mi_F & \mi_F & \mi_F\\
\mi_F & \mi_F & \mi_F & \mi_F\\ \hline
\mi_F^2 & \mi_F^2 & \mi_F & \mi_F\\
\mi_F^2 & \mi_F^2 & \mi_F & \mi_F\\
\end{array}
\right).
\end{eqnarray}

Put
$n = 8m-4$ and $\beta = \beta(0)$.
\begin{lem}\label{lem:semisimple}
Suppose that an irreducible smooth representation
$\pi$ of $G$ contains a skew stratum $[\Lambda', n', n'-1, \beta(c)]$, for some $c \in \ri_0$.
Then $\pi$ contains $[\Lambda, n, n-1, \beta]$.
\end{lem}
\begin{proof}
The assertion follows from the inclusion
$\beta(c) + \rad_{1-n'}(\Lambda') \subset
\beta + \rad_{1-n}(\Lambda)$.
\end{proof}

\begin{rem}
Let $\pi$ be an irreducible smooth representation of $G$.
By Lemma~\ref{lem:semisimple} and 
the argument in \S \ref{subsec:norm},
$\pi$ is of level $n/2$ and of characteristic polynomial 
$(X -\alpha)^4 \pmod{\mi_F}$
if and only if $\pi$ contains $[\Lambda, n, n-1, \beta]$.
\end{rem}

We are going to classify the irreducible smooth representations
of $G$
of level $n/2$ and with characteristic polynomial 
$(X -\alpha)^4 \pmod{\mi_F}$.
Let $E = F[\beta]$ and 
let $B$ denote
the $A$-centralizer of $\beta$.
Then $E$ and $B$ are $\sigma$-stable,
and the algebra $E$ is 
a totally ramified extension of degree 2 over $F$
with uniformizer $\p_E = \p^m \beta \e$.
Put $E_0 =\{ x\in E\ |\ \sigma(x) = x \}$ and 
$G' = G\cap B$.
Then $E$ is the quadratic unramified extension over $E_0$
and 
$G'$ is the unramified unitary group over $E_0$
corresponding to the involutive algebra $(B, \sigma)$.
Note that $\p_E \in E_0$.

Observe that 
$\beta \Lambda(i) = \Lambda(i-n)$, $i \in \Z$.
Then 
$\Lambda$ is an $\ri_E$-lattice sequence in $V$
of $\ri_E$-period 4
since $\ri_E = \ri_F[\p_E]$.
Since $\beta^2 \in F^\times$, the map
$s: A \rightarrow B$ defined by
\begin{eqnarray}
s(X) = (X + \beta X \beta^{-1})/2,\ X \in A
\end{eqnarray}
is a $(B, B)$-bimodule projection.
Let $B^\bot$ denote the orthogonal complement of $B$
in $A$ with respect to the pairing induced by
$\mathrm{tr}_{A/F}$.
We have $B^\bot = \ker s$ and 
$A = B \oplus B^\bot$.
Note that the set $B^\bot$
is also $\sigma$-stable.

For 
$k \in \Z$,
we abbreviate
$\rad_k = \rad_k(\Lambda)$,
$\rad'_k =  \rad_k \cap B$ and
$\rad_k^\bot = \rad_k(\Lambda) \cap B^\bot$.
\begin{prop}\label{prop:decomp_E}
For all $k \in \Z$, we have
$\rad_k = \rad'_k \oplus 
\rad_k^\bot$.
\end{prop}
\begin{proof}
Since $\beta$ normalizes $\rad_k$,
we have $s(\rad_k) = \rad'_k$.
The proposition follows from the equation
$X = s(X) +(X-s(X))$, $X \in A$.
\end{proof}

As in \cite{GSp4} (4.16),
we define a $\sigma$-stable $\ri_F$-lattice in $A$ by
\begin{eqnarray}\label{eq:def_J}
\J = \rad'_n\oplus \rad_{[n/2]+1}^\bot
\end{eqnarray}
and an open compact subgroup $J = (1 + \J)\cap G$ of $G$.
Then we get
$P_{\Lambda, n} \subset J \subset P_{\Lambda, [n/2]+1}$
and hence the quotient $J/P_{\Lambda, n+1}$ is abelian.
As usual, there is an isomorphism 
\begin{eqnarray*}
\g_{\Lambda, -n}/ \g\cap \J^* \simeq (J/P_{\Lambda, n+1})^\wedge;\
b +\g \cap \J^* \mapsto \Psi_b,
\end{eqnarray*}
given by
\begin{eqnarray}\label{eq:Psi}
\Psi_b(x) = \Omega(\mathrm{tr}_{A/F_0}(b(x-1))),\ x \in J.
\end{eqnarray}
By Proposition~\ref{prop:decomp_E} and \cite{BK2} (2.10),
we have
\begin{eqnarray}
\J^* = \rad_{1-n}' \oplus \rad_{-[n/2]}^\bot.
\end{eqnarray}

\begin{rem}
In general, we need to consider the  lattice
$\rad'_n\oplus \rad_{[(n+1)/2]}^\bot$.
But in this case,
it follows from (\ref{eq:E_fil}) that
$\rad_{[n/2]+1}^\bot = 
\rad_{[(n+1)/2]}^\bot$.
\end{rem}

We also abbreviate
$\g_k = \rad_k \cap \g$,
$\g'_k = \rad'_k \cap \g$ and
$\g_k^\bot = \rad_k^\bot \cap \g$, for $k \in \Z$.
For $X \in A$,
we write $\mathrm{ad}(\beta)(X) = \beta X -X\beta$. 
Since $\beta \in \g'_{-n}$,
the map $\ad(\beta)$ induces a quotient map
\[
\ad(\beta):
\g_{k}^\bot/\g_{k+1}^\bot \rightarrow 
\g_{k-n}^\bot /\g_{k-n+1}^\bot,\ k \in \Z.
\]
\begin{lem}\label{lem:u11}
For $k \in \Z$,
the map
$\ad(\beta):
\g_{k}^\bot/\g_{k+1}^\bot \simeq 
\g_{k-n}^\bot /\g_{k-n+1}^\bot$ is an isomorphism.
\end{lem}
\begin{proof}
By the periodicity of $\{\seqrad_k\}_{k\in \Z}$,
it suffices to prove that
$\mathrm{ad}(\beta)$ induces an injection, for all $k \in \Z$.
Let
$X \in \g_k^\bot$ satisfy
$\mathrm{ad}(\beta)(X) \in \g_{k-n+1}^\bot$.
Since $\beta^{-1} \in \rad_{n}$,
we have
$X = X-s(X) = -\ad(\beta)(X)\beta^{-1}/2
\in \g_{k+1}^\bot$,
 as required.
\end{proof}

\begin{prop}[\cite{GSp4} Lemma 4.4]\label{prop:5c1_surj}
Suppose that an element 
$\gamma \in \beta + \g_{\Lambda, 1-n}$ 
lies in $\g\cap B$ 
modulo $\g_{\Lambda, k-n}$
for some integer $k \geq 1$.
Then, 
there exists 
$p \in P_{\Lambda, k}$ such that
$\mathrm{Ad}(p)(\gamma)  \in \g\cap B$.
\end{prop}
\begin{proof}
Exactly the same as the proof of \cite{GSp4} Lemma 4.4.
\end{proof}

As an immediate corollary of the proof we have
\begin{cor}[\cite{GSp4} Corollary 4.5]\label{cor:5c1_J}
$\Ad(J)(\beta + \g'_{1-n}) = \beta + \g \cap \J^*$.
\end{cor}

\begin{prop}[\cite{GSp4} Theorem 4.1]\label{prop:5_2}
Let $\pi$ be an irreducible smooth representation of $G$.
Then $\pi$ contains $[\Lambda, n, n-1, \beta]$ if and only if
the restriction of $\pi$ to $J$ contains $\Psi_\beta$.
\end{prop}
\begin{proof}
Since $\Psi_\beta$ is an extension of $\psi_\beta$ to $J$,
$\pi$ contains $[\Lambda, n, n-1, \beta]$ if 
$\pi$ contains $\Psi_\beta$.

Suppose that $\pi$ contains $[\Lambda, n, n-1, \beta]$.
Then we can find
$\gamma \in \beta + \g_{\Lambda, 1-n}$
so that
a skew stratum $[\Lambda, n, [n/2], \gamma]$ 
occurs in $\pi$.
By Proposition~\ref{prop:5c1_surj},
after conjugation by some element in $P_{\Lambda, 1}$,
we may assume 
$\gamma$ lies in $\beta + \g'_{1-n}$.
This implies that
the restriction of $\psi_\gamma$ to $J$
is $\Psi_\beta$.
This completes the proof.
\end{proof}

We write $\Psi = \Psi_\beta$,
$J' = J\cap G'$ and $\Psi' = \Psi|_{J'}$.
\begin{thm}\label{thm:E-isom}
With the notation as above,
there is a  $*$-isomorphism 
$\eta : \He(G'//J', \Psi') \simeq \He(G//J, \Psi)$
which preserves support,
that is, $\eta$ satisfies
$\mathrm{supp}(\eta(f)) = J\mathrm{supp}(f)J$, for $f \in \He(G'//J', \Psi')$.
\end{thm}
\begin{rem}\label{rem:u11_E}
The isomorphism $\eta$ allows us to identify  
$\Irr(G')^{(J', \Psi')}$ with $\Irr(G)^{(J, \Psi)}$.
Since $\Psi'$ can extend to a character of $G'$,
we obtain a  bijection from $\Irr(G')^{(J', 1)}$
to $\Irr(G)^{(J, \Psi)}$,
where $1$ denotes the trivial representation of $J'$.
Since the centers of  $G$ and $G'$ are compact,
the support preservation of $\eta$ implies
that 
this map preserves the supercuspidality of representations.
\end{rem}

We commence the proof of Theorem~\ref{thm:E-isom}.
\begin{prop}\label{prop:E_surj}
$I_G(\Psi) = JG' J$.
\end{prop}
\begin{proof}
Note that an element $g$ in $G$ lies in 
$I_G(\Psi)$ if and only if
$\mathrm{Ad}(g)(\beta +\g \cap \J^*)\cap 
(\beta +\g \cap \J^*) \neq \emptyset$ .
Clearly, we have
$J I_G(\Psi) J = I_G(\Psi)$,
$G' \subset I_G(\Psi)$ and hence
$J G' J \subset I_G(\Psi)$.

Let $g \in I_G(\Psi)$.
Then Corollary~\ref{cor:5c1_J} says that
there exists an element $k$ in $JgJ$
such that
$\Ad(k)(\beta+\g'_{1-n})\cap (\beta +
\g'_{1-n}) \neq \emptyset$.
Then there are $x, y$ in $\g'_{1-n}$
such that
$\ad(\beta)(k) = kx-yk$.
Projecting this equation on $B^\bot$,
we have $\ad(\beta)(k^\bot) = k^\bot x-yk^\bot$,
where $k^\bot$ denote the $B^\bot$-component
of $k$.
Suppose $k^\bot \in \rad_{l}$, for some $l \in \Z$.
Then we have
$\ad(\beta)(k^\bot) \in \rad_{l-n+1}$,
and by  applying the proof of Lemma~\ref{lem:u11},
$k^\bot \in \rad_{l+1}$.
So we conclude $k^\bot = 0$
and hence $k \in G'$.
This completes the proof.
\end{proof}

For $x$ in $A$,
we denote by $x'$ its $B$-component
and by $x^\bot$ its $B^\bot$-component.
\begin{lem}\label{lem:intertwine_J}
For $g \in JG'J$,
we have
$\nu_\Lambda(g^\bot) \geq \nu_\Lambda(g') +[n/2]+1$.
\end{lem}
\begin{proof}
Put $k = \nu_\Lambda(g)$.
Then, for any element $y$ in $JgJ$,
we have
$y \equiv g \pmod{\rad_{k+[n/2]+1}}$,
so that
$y \equiv g^\bot \pmod{B+ \rad_{k +[n/2]+1}}$.
Therefore if $g \in JG'J$,
then $g^\bot \in \rad_{k +[n/2]+1}^\bot$.
Now the lemma follows immediately.
\end{proof}

\begin{prop}\label{prop:E_inj}
For
$g \in G'$,
we have
$JgJ \cap G' = J' g J'$.
\end{prop}
\begin{proof}
Since $\beta\sigma(\beta) = -\beta^2 \in F^\times$,  
the map $\Ad(\beta)$
gives an automorphism on $G$ of order 2
and $J$ is a pro-$p$ subgroup of $G$ closed under $\Ad(\beta)$.
The assertion follows from \cite{St2} Lemma 2.1.
\end{proof}
We can identify 
$M_2(E)$ with $B$ as follows:
For $a, b, c, d, x, y, z, w \in F$,
\begin{eqnarray*}
X = \left(
\begin{array}{cc}
a + \p_E x & b +\p_E y\\
c +\p_E z & d +\p_E w
\end{array}
\right)
\mapsto
\left(
\begin{array}{cc|cc}
a & -b\e^{-1} & x\e & y\\
-c \e  & d & -\varepsilon z & -w\e \\ \hline
\p \alpha x \e  & -\p \alpha y & a & b\e^{-1}\\
\p \alpha \varepsilon z & -\p \alpha w \e  & c \e  & d
\end{array}
\right).
\end{eqnarray*}
The involution $\sigma$ on $M_2(E)$ induced by 
this identification
maps $X$ to 
\begin{eqnarray*}
\sigma(X) = \left(
\begin{array}{cc}
\overline{d} + \p_E\overline{w} & \overline{b} + \p_E\overline{y}\\
\overline{c} + \p_E\overline{z} & \overline{a} + \p_E\overline{x}
\end{array}
\right).
\end{eqnarray*}
So the group
 $G'$ is isomorphic to the unramified $U(1,1)$
over $E_0$.

We write $B' = G'\cap P_{\Lambda, 0}$.
Since the $\ri_E$-lattice sequence $\Lambda$ satisfies
$\Lambda(0) \supsetneq  \Lambda(3) 
\supsetneq \p_E\Lambda(0)$,
the group
$B'$ is the Iwahori subgroup of $G'$
and 
$J' 
= G' \cap (1+\p_E^{2m-1} \rad_0(\Lambda))$
is its ($4m-2$)-th standard congruence subgroup.
Set
\begin{eqnarray*}
s_1 = \left(
\begin{array}{cc}
0 & 1\\
1 & 0
\end{array}
\right),\
s_2 = \left(
\begin{array}{cc}
0 & \p_E^{-1}\\
\p_E & 0
\end{array}
\right) \in G',
\end{eqnarray*}
$S = \{s_1, s_2\}$,
and $W' = \langle S \rangle$.
Then we have a Bruhat decomposition 
$G' = B' W' B'$.
\begin{lem}\label{lem:cosets}
Let $t$ be a non negative integer.
Then we have

(i)
$[J' (s_1s_2)^t J': J'] = 
[J'(s_2 s_1)^tJ' : J'] 
= q^{2t}$,
$[J' (s_1s_2)^ts_1 J': J'] =  [J' (s_2 s_1)^t s_2 J' : J'] = q^{2t+1}$.

(ii)
$[J (s_1s_2)^t J: J] = [J(s_2 s_1)^tJ : J]  = q^{4t},$
$[J (s_1s_2)^ts_1 J: J] = q^{4t+1}$,
$[J (s_2 s_1)^t s_2 J: J] = q^{4t+3}$.
\end{lem}
\begin{proof}
For $g$ in $G'$,
we have
$[J'gJ' : J'] = [J':J'\cap g J' g^{-1}] =
[\g'_n: \g_n'\cap g\g_n'g^{-1}]$
and
\begin{eqnarray*}
[JgJ : J] & = & [J:J\cap gJg^{-1}] = 
[\g \cap \J : \g \cap \J\cap g\J g^{-1}]\\
& = & 
[\g_n': \g_n'\cap g\g_n'g^{-1}]
\cdot
[\g_{[n/2]+1}^\bot:
\g_{[n/2]+1}^\bot \cap g\g_{[n/2]+1}^\bot g^{-1}].
\end{eqnarray*}
These can be
directly calculated
for $w \in W'$ by (\ref{eq:E_fil}). 
\end{proof}

For $\mu \in E$ and $\nu \in E^\times$, we 
define 
\begin{eqnarray*}
u(\mu) = 
\left(
\begin{array}{cc}
1 & \mu\\
0 & 1
\end{array}
\right),\ 
\underline{u}(\mu) = 
\left(
\begin{array}{cc}
1 & 0\\
\mu & 1
\end{array}
\right),\ \mathrm{and}\
h(\nu) = 
\left(
\begin{array}{cc}
\nu & 0\\
0 & \sigma(\nu)^{-1}
\end{array}
\right).
\end{eqnarray*}
For $g \in G'$,
let $e_g$ denote the element in 
$\He(G'//J', \Psi')$
such that $e_g(g) = 1$ and $\mathrm{supp}(e_g) = J'gJ'$.
\begin{thm}\label{thm:e}
The algebra $\He(G'//J', \Psi')$ is generated by the elements $e_g$, $g \in B' \cup S$.
These elements satisfy the following relations:

(i) $e_k = \Psi(k) e_1$, $k \in J'$,

(ii) $e_k* e_{k'} = e_{kk'}$, $k, k' \in B'$,

(iii) $e_s * e_k = e_{sks}* e_s$, $s \in S$, $k \in sB's\cap B'$,

(iv) 
$e_{s_1}* e_{s_1}  =   [J's_1 J': J'] \sum_{x \in \ri_{E_0}/\mi_{E_0}} 
e_{\underline{u}(\varpi_E^{2m-1} x \sqrt{\varepsilon})}$,

$e_{s_2}* e_{s_2}  =  [J's_2 J': J'] \sum_{x \in \ri_{E_0}/\mi_{E_0}} 
e_{{u}(\varpi_E^{2m-2} x \e)}$,

(v)
For $\mu \in \ri_{E_0}^\times \e$, 
\begin{eqnarray*}
e_{s_1}* e_{u(\mu)} * e_{s_1} & = & [J's_1J':J']
e_{u(\mu^{-1})}* e_{s_1}* e_{h(\mu)}* e_{u(\mu^{-1})},\\
e_{s_2}* e_{\underline{u}(\varpi_E \mu)} * e_{s_2} & = & [J's_2J':J']
e_{\underline{u}(\varpi_E \mu^{-1})}* e_{s_2}* e_{h(-\mu^{-1})}* 
e_{\underline{u}(\varpi_E \mu^{-1})}.
\end{eqnarray*}
The above relations are a defining set for this algebra.
\end{thm}

\begin{proof}
The proof is very similar to the proof of 
\cite{Harish} Chapter 3 Theorem 2.1.
\end{proof}

For $g \in G'$,
let $f_g$ denote the element in $\He(G//J, \Psi)$ such that
$f_g(g) = 1$ and $\mathrm{supp}(f_g) = JgJ$.
\begin{thm}\label{thm:f}
The algebra $\He(G//J, \Psi)$ is generated by
 $f_g$, $g \in S \cup B'$ and satisfies the following relations:

(i) $f_k = \Psi(k) f_1$, $k \in J'$,

(ii) $f_k* f_{k'} = f_{kk'}$, $k, k' \in B'$,

(iii) $f_s * f_k = f_{sks}* f_s$, $s \in S$, $k \in sB's\cap B'$,

(iv) 
$f_{s_1}* f_{s_1}  =   [Js_1 J: J] \sum_{x \in \ri_{E_0}/\mi_{E_0}} 
f_{\underline{u}(\varpi_E^{2m-1} x \sqrt{\varepsilon})}$,

$f_{s_2}* f_{s_2}  =  [Js_2 J: J] \sum_{x \in \ri_{E_0}/\mi_{E_0}} 
f_{{u}(\varpi_E^{2m-2} x \e)}$,

(v)
For $\mu \in \ri_{E_0}^\times \e$, 
\begin{eqnarray*}
f_{s_1}* f_{u(\mu)} * f_{s_1} & = & [Js_1J:J]
f_{u(\mu^{-1})}* f_{s_1}* f_{h(\mu)}* f_{u(\mu^{-1})},\\
f_{s_2}* f_{\underline{u}(\varpi_E \mu)} * f_{s_2} & = & q^2
f_{\underline{u}(\varpi_E \mu^{-1})}* f_{s_2}* f_{h(-\mu^{-1})}* 
f_{\underline{u}(\varpi_E \mu^{-1})}.
\end{eqnarray*}
\end{thm}
\begin{proof}
Recall that if $x,  y \in G'$ satisfy $[JxJ:J] 
 [JyJ:J] =[JxyJ:J]$,
 then
$JxJyJ = JxyJ$ and
$f_x* f_y = f_{xy}$.

By Proposition~\ref{prop:E_surj}, 
$\He(G//J, \Psi)$ is linearly spanned by 
$f_g$, $g \in G'$.
For $g \in G'$,
we write $g = b_1 w b_2$ where $b_1, b_2 \in B'$ and
$w \in W'$.
Then we have $f_g = f_{b_1}*f_w*f_{b_2}$ because
$B'$ normalizes $J$.
Let $w = s_{i_1} s_{i_2}\cdots s_{i_l}$ be a minimal expression
for $w$
with $s_{i_j} \in S$.
It follows from Lemma~\ref{lem:cosets}
that $f_w = f_{s_{i_1}}*f_{s_{i_2}}*\cdots *f_{s_{i_l}}$.
Therefore 
$\He(G//J, \Psi)$ is generated by $f_g$, $g \in B' \cup S$.

 Relations (i), (ii), and (iii) are obvious.
Since 
$[Js_1J:J] = [J's_1J':J']$,
we can take a common system of representatives for
$J/J\cap sJs$ and $J'/J'\cap sJ's$.
Then the proof of relations on $s_1$ is also obvious.

We shall give the proof of 
 relations on $s_2$ only in the case when $m = 2k+1$, $k \geq 0$,
 and omit the easy modification when $m$ is even.

We will abbreviate $s = s_2$.
By the identification $M_2(E) = B$, we have
\begin{eqnarray}
s
= \left(
\begin{array}{cccc}
0 & 0 &0 & \p^{-1}\alpha^{-1}\varepsilon^{-1}\\
0 & 0 &  -\varepsilon& 0\\
0 & -\varepsilon^{-1} & 0 &0\\
\p \alpha \varepsilon & 0 & 0 & 0
\end{array}
\right).
\end{eqnarray}
We can take a  common
system of representatives for $J/J\cap sJs$
and $J\cap sJs\backslash J$ as
\begin{eqnarray}
x(a, b, c)
= \left(
\begin{array}{cccc}
1 & 0 &0 & 0\\
\p^{k+1}a - \p^m c\varepsilon & 1 &  \p^k \alpha^{-1}b& 0\\
0 & 0 & 1 &0\\
\p^{k+1} b & 0 & \p^{k+1}a +\p^m c\varepsilon & 1
\end{array}
\right),\ 
\end{eqnarray}
where
$a, b \in \ri_0\e/\mi_0\e,\ c \in \ri_0/\mi_0.$
Note that $x(a,b,c)$ lies in the kernel of $\Psi$.

(iv)
For $g \in G$,
let $\delta_g$ denote the unit point mass at $g$.
Then we have
\[
f_s = \sum_{x \in J\cap s Js\backslash J} 
\Psi^{-1}(x)
f_1*\delta_{sx}
=
\sum_{
y \in J/J\cap sJs} 
\Psi^{-1}(y)
\delta_{ys}*f_1,
\]
so that
\begin{eqnarray*}
f_s * f_s
& = & \sum_{x \in J\cap s Js\backslash J,\
y \in J/J\cap sJs} 
\Psi^{-1}(xy)
f_1*\delta_{sxys}*f_1\\
 & = & [JsJ:J] 
 \sum_{x \in J\cap s Js\backslash J} \Psi^{-1}(x) f_1*\delta_{sxs}*f_1.
\end{eqnarray*}

We remark that for $g \in G$,
\begin{eqnarray*}
f_1*\delta_g *f_1
= 
\left\{
\begin{array}{cl}
[JgJ:J]^{-1} f_g, & \mathrm{if}\ g \in I_G(\Psi),\\
0, & \mathrm{if}\ g \not\in I_G(\Psi).
\end{array}
\right.
\end{eqnarray*}
The $B$-component of $sx(a,b,c)s$ 
is $sx(0,0,c)s$ and lies in $\rad_0$.
Therefore,
if $sx(a,b,c)s \in I_G(\Psi) = JG' J$,
then
by Lemma~\ref{lem:intertwine_J},
$sx(a,b,0)s \in 1 +\rad_{[n/2]+1}$
and hence $a \equiv b \equiv 0$.

So we obtain 
\begin{eqnarray*}
f_s * f_s
 & = & 
 [JsJ:J] 
 \sum_{c\in \ri_0/\mi_0} 
 f_1*\delta_{sx(0, 0, c)s}*f_1\\
& = &
 [Js J :J]\sum_{c \in \ri_0/\mi_0} 
f_1*\delta_{u(\p_E^{-2}\p^m c\e)}*f_1\\
& = &
[Js J :J]\sum_{x \in \ri_{E_0}/\mi_{E_0}} 
f_{u(\p_E^{2m-2}x\e)}.
\end{eqnarray*}

(v)
Let $\mu \in \ri_{E_0}^\times \e$.
Put 
$u = \underline{u}(\varpi_E \mu) \in B'$.
Then
$f_u = f_1*\delta_u = \delta_u *f_1$.
Since
$u$ normalizes the pair $(J, \Psi)$,
we have
\begin{eqnarray*}
f_s * f_u*f_s
& = & \sum_{x \in J\cap s Js\backslash J,\
y \in J/J\cap sJs} 
\Psi^{-1}(xy)
f_1*\delta_{sxuys}*f_1\\
 & = & [JsJ:J] 
 \sum_{x \in J\cap s Js\backslash J} \Psi^{-1}(x) f_1*\delta_{sxus}*f_1.
\end{eqnarray*}

The $B$-component of $sx(a,b,c)us$, which is
identified with $u(\p_E^{-1}\mu + \p^m \p_E^{-2}c\e)$,
lies in $\rad_{-3}$.
Suppose $sx(a,b,c)us \in I_G(\Psi) = JG' J$.
Then
by Lemma~\ref{lem:intertwine_J},
we have $sx(a,b,0)s \in 1+\rad_{8k}$
and hence
$b \equiv 0$.
So we obtain 
\begin{eqnarray*}
f_s * f_u*f_s
 & = & [JsJ:J] 
 \sum_{a \in \ri_0\e/\mi_0\e,\ c\in \ri_0/\mi_0} 
 f_1*\delta_{sx(a, 0, c)us}*f_1.
\end{eqnarray*}

For the moment,
we fix $a$ and $c$.
We put $\nu = \mu +\p^m \p_E^{-1}c\e$,
$v = \underline{u}(\p_E \nu^{-1})$,
$h = h(-\nu^{-1})$, and $x= x(a,0,0)$.
Then $sx(a,0,c)us = sxs s \underline{u}(\p_E \nu)s
= sxs vshv
= [sxs, v] vshv (hv)^{-1}xhv$.
Since $hv \in B'$,
we have $(hv)^{-1}xhv \in J$
and $\Psi((hv)^{-1}xhv) = \Psi(x) = 1$.
Since $sxs \in P_{\Lambda, 8k+1}$
and $v \in P_{\Lambda, 3}$,
we obtain 
$[sxs, v] \in P_{\Lambda, [n/2]+1}$.
If we write $sxs = 1+y$ and $v = 1 +z$,
then we have
$[sxs, v] = (1+y)(1+z)(1-y)(1-z) 
\equiv 1 +yzy + yzyz \pmod{B^\bot}$.
So the $B$-component of $[sxs, v]$ lies in
$1 +\rad_{16k+5} = 1+\rad_{n+1}$.
This implies that $[sxs, v] \in J$ and 
$\Psi([sxs,v]) =1$.

We therefore have
\begin{eqnarray*}
 f_1*\delta_{sx(a, 0, c)us}*f_1
& = & f_1 * \delta_{vshv}*f_1\\
& = & [JsJ:J]^{-1} f_v*f_s*f_h*f_v\\
& = & 
[JsJ:J]^{-1}
f_{\underline{u}(\varpi_E \mu^{-1})}* f_{s}* f_{h(-\mu^{-1})}* 
f_{\underline{u}(\varpi_E \mu^{-1})},
\end{eqnarray*}
for any $a \in \ri_0\e/\mi_0\e$ and $c \in \ri_0/\mi_0$.
So we conclude that
\begin{eqnarray*}
 f_s*f_u*f_s
& = & 
q^2
f_{\underline{u}(\varpi_E \mu^{-1})}* f_{s}* f_{h(-\mu^{-1})}* 
f_{\underline{u}(\varpi_E \mu^{-1})},
\end{eqnarray*}
as required.
\end{proof}

We return to the proof of Theorem~\ref{thm:E-isom}.
By Theorems~\ref{thm:e} and \ref{thm:f},
there is an algebra homomorphism
$\eta: \He(G'//J', \Psi) \rightarrow
\He(G//J, \Psi)$ induced by
\[
\eta(e_{s_1}) = f_{s_1},\
\eta(e_{s_2}) = q^{-1} f_{s_2},\
\eta(e_b) = f_b,\ b \in B'.
\]
Note that,
for $g \in G'$, 
$\eta$ maps $e_g$ to $(\mathrm{vol}(J'gJ')/\mathrm{vol}(JgJ))^{1/2}f_g$.
Then, it follows from Proposition~\ref{prop:E_surj} that
$\eta$ is surjective.
Proposition~\ref{prop:E_inj} implies that
$\eta$ is injective.
The $*$-preservation of $\eta$ is obvious
since
$e_g^* = e_{g^{-1}}$ and $f_g^* = f_{g^{-1}}$,
for $g \in G'$.

\subsection{Case (\ref{sec:h_int}d)}\label{5_1}
Let 
$\Lambda'$ denote the $\ri_F$-lattice sequence 
in (\ref{eq:st_20}).
We fix an element $d \in \ri_0^\times$
and a positive odd integer $n' = 2m-1$.
For $b, c \in \ri_0$,
we set
\begin{eqnarray}\label{eq:beta1}
\beta(b, c) = \p^{-m}
\left(
\begin{array}{cccc}
0 & 0 & 0 & \e\\
0 & 0 & b\e & 0\\ 
0 & \p c\e & 0 & 0\\
\p d\e & 0 & 0 & 0
\end{array}
\right) \in \g_{\Lambda', -n'}.
\end{eqnarray}
In this section,
we consider the irreducible smooth representations of $G$
which contain a skew stratum $[\Lambda', n', n'-1, \beta]$,
where $\beta$ equals to 
$\beta(0, 0)$, $\beta(0, 1)$, or $\beta(1,0)$.

We define a self-dual $\ri_F$-lattice sequence $\Lambda$ 
in $V$ with $e(\Lambda) = 8$ and $d(\Lambda) = 3$
by
\[
\Lambda(0) = \Lambda(1) = \Lambda(2) = \Lambda(3) = N_0,\
\Lambda(4) = N_1,\
\Lambda(5) = \Lambda(6) = N_2,\
\Lambda(7) = \p N_1^\#,
\]
\begin{eqnarray}\label{eq:L1}
\Lambda(i +8k) = \p^k \Lambda(i),\ 0 \leq i \leq 7,\ 
k \in \Z.
\end{eqnarray}
The filtration
 $\{\rad_k(\Lambda)\}_{k \in \Z}$ 
is given by the following:
\[
\rad_0(\Lambda) = 
\left(
\begin{array}{c|cc|c}
\ri_F & \ri_F & \ri_F & \ri_F\\ \hline
\mi_F & \ri_F & \ri_F & \ri_F\\ 
\mi_F & \mi_F & \ri_F & \ri_F\\ \hline
\mi_F & \mi_F & \mi_F & \ri_F
\end{array}
\right),\
\rad_1(\Lambda) = 
\left(
\begin{array}{c|cc|c}
\mi_F & \ri_F & \ri_F & \ri_F\\ \hline
\mi_F & \mi_F & \ri_F & \ri_F\\ 
\mi_F & \mi_F & \mi_F & \ri_F\\ \hline
\mi_F & \mi_F & \mi_F & \mi_F
\end{array}
\right),
\]
\[
\rad_2(\Lambda) = 
\left(
\begin{array}{c|cc|c}
\mi_F & \mi_F & \ri_F & \ri_F\\ \hline
\mi_F & \mi_F & \ri_F & \ri_F\\ 
\mi_F & \mi_F & \mi_F & \mi_F\\ \hline
\mi_F & \mi_F & \mi_F & \mi_F
\end{array}
\right),\
\rad_3(\Lambda) = 
\left(
\begin{array}{c|cc|c}
\mi_F & \mi_F & \ri_F & \ri_F\\ \hline
\mi_F & \mi_F & \mi_F & \ri_F\\ 
\mi_F & \mi_F & \mi_F & \mi_F\\ \hline
\mi_F & \mi_F & \mi_F & \mi_F
\end{array}
\right),
\]
\[
\rad_4(\Lambda) = 
\left(
\begin{array}{c|cc|c}
\mi_F & \mi_F & \mi_F & \ri_F\\ \hline
\mi_F & \mi_F & \mi_F & \mi_F\\ 
\mi_F & \mi_F & \mi_F & \mi_F\\ \hline
\mi_F & \mi_F & \mi_F & \mi_F
\end{array}
\right),\
\rad_5(\Lambda) = 
\left(
\begin{array}{c|cc|c}
\mi_F & \mi_F & \mi_F & \mi_F\\ \hline
\mi_F & \mi_F & \mi_F & \mi_F\\ 
\mi_F & \mi_F & \mi_F & \mi_F\\ \hline
\mi_F^2 & \mi_F & \mi_F & \mi_F
\end{array}
\right),
\]
\begin{eqnarray}\label{eq:filt-5c1}
& \rad_6(\Lambda) = 
\left(
\begin{array}{c|cc|c}
\mi_F & \mi_F & \mi_F & \mi_F\\ \hline
\mi_F & \mi_F & \mi_F & \mi_F\\ 
\mi_F^2 & \mi_F & \mi_F & \mi_F\\ \hline
\mi_F^2 & \mi_F^2 & \mi_F & \mi_F
\end{array}
\right),\
\rad_7(\Lambda) = 
\left(
\begin{array}{c|cc|c}
\mi_F & \mi_F & \mi_F & \mi_F\\ \hline
\mi_F & \mi_F & \mi_F & \mi_F\\ 
\mi_F^2 & \mi_F^2 & \mi_F & \mi_F\\ \hline
\mi_F^2 & \mi_F^2 & \mi_F & \mi_F
\end{array}
\right).
\end{eqnarray}

If we write
$n = 8m-4$ and $\beta = \beta(0,0)$,
then we obtain the analogue of Lemma~\ref{lem:semisimple},
that is,
\begin{eqnarray}\label{eq:similitude}
\Irr(G)^{(P_{\Lambda', n'}, \psi_{\beta(0,1)})}
\cup \Irr(G)^{(P_{\Lambda', n'}, \psi_{\beta})}
\subset \Irr(G)^{(P_{\Lambda, n}, \psi_{\beta})}.
\end{eqnarray}

Take $\lambda \in \ri_F$ such that
$\lambda\overline{\lambda} = d$
and put
\[
t = 
\left(
\begin{array}{cccc}
0 & 0 & 0 & \lambda^{-1}\\
0 & 0 & 1 & 0\\
0 & \p & 0 & 0\\
\p \overline{\lambda} & 0 & 0 & 0
\end{array}
\right).
\]
Then $t$ is a similitude on $(V, f)$ and hence
$t$ acts on the set of 
skew strata as well.
Note that $t\Lambda'$ is a translate of $\Lambda'$.
As a conjugate of (\ref{eq:similitude}) by $t$,
we get 
\begin{eqnarray}
\Irr(G)^{(P_{\Lambda', n'}, \psi_{\beta(1,0)})}
\cup \Irr(G)^{(P_{\Lambda', n'}, \psi_{\beta})}
\subset \Irr(G)^{(P_{t\Lambda, n}, \psi_{\beta})}.
\end{eqnarray}
Since $\Irr(G)^{(P_{t\Lambda, n}, \psi_{\beta})}$
is the $t$-conjugate of $\Irr(G)^{(P_{\Lambda, n}, \psi_{\beta})}$,
there is a bijection  
$\Irr(G)^{(P_{\Lambda, n}, \psi_{\beta})} \simeq
\Irr(G)^{(P_{t\Lambda, n}, \psi_{\beta})}$;
 $\pi \mapsto \pi^t$.
We shall concentrate on classifying
the irreducible smooth representations of
$G$ containing 
$[\Lambda, n, n-1, \beta]$.
\begin{rem}
We note that 
$\Irr(G)^{(P_{t\Lambda, n}, \psi_{\beta})}
\cap \Irr(G)^{(P_{\Lambda, n}, \psi_{\beta})}
\supset \Irr(G)^{(P_{\Lambda', n'}, \psi_{\beta})}$.
\end{rem}

Define an orthogonal $F$-splitting $V = V^1 \bot V^2$ by
$V^1 = Fe_2 \oplus Fe_3$ and 
$V^2 = Fe_1 \oplus Fe_4$.
Then
the skew stratum $[\Lambda, n, n-1, \beta]$
is split with respect to $V = V^1 \bot V^2$.
Put $\beta_j = \beta|_{V^j}$, for $j = 1,2$.
Then $\beta_1 = 0$ and
$E_2 = F[\beta_2]$ is a totally ramified extension 
of degree 2 over $F$.
We denote by 
$B$ the $A$-centralizer of $\beta$
and by $B^\bot$ its orthogonal complement 
with respect to 
the pairing induced by
$\mathrm{tr}_{A/F}$.
Then $B$ and $B^\bot$ are $\sigma$-stable.
For $i, j \in \{ 1, 2 \}$,
we write $A^{ij} = \mathrm{Hom}_F(V^j, V^i)$.
Since $\beta_2^2 \in F^\times \cdot 1_{V^2}$,
the map 
\[
s_2: A^{22} \rightarrow  E_2 ;
X \mapsto (X + \beta_2 X \beta_2^{-1})/2,\
X \in A^{22}
\]
is an $(E_2, E_2)$-bimodule projection.
So we get
\begin{eqnarray}\label{eq:B_bot}
B = A^{11}\oplus E_2,\
B^\bot = 
A^{12} \oplus A^{21} \oplus
\ker s_2,\
A = B \oplus B^\bot.
\end{eqnarray}

For brevity,
we write
$\rad_k = \rad_k(\Lambda)$,
$\rad_k' = \rad_k\cap B$, 
$\rad_k^\bot = \rad_k\cap B^\bot$,
$\g'_k = \g_{\Lambda, k}\cap B$
and 
$\g_k^\bot = \g_{\Lambda, k}\cap B^\bot$,
for $k \in \Z$.
\begin{prop}\label{prop:5c1-eq}
For $k \in \Z$, we have
$\rad_k = \rad_k' \oplus \rad_k^\bot$.
\end{prop}
\begin{proof}
By \cite{BK2} (2.9),
we have
$\seqrad_k = \bigoplus_{1 \leq i,j\leq 2} 
\seqrad_k \cap A^{ij}$
and $\rad_k\cap A^{ii} = \rad_k(\Lambda^i)$, for $i = 1, 2$.
Since $\beta_2$ normalizes 
$\rad_k(\Lambda^2)$,
we have
$\rad_k(\Lambda^2) = (\seqrad_k \cap E_2)
\oplus (\seqrad_k \cap \ker s_2)$.
Now the proposition follows from
(\ref{eq:B_bot}).
\end{proof}
\begin{rem}
Since
$[(n+1)/2] \equiv 2$ or $6 \pmod{8}$,
it follows from (\ref{eq:filt-5c1}) that
$\rad_{[(n+1)/2]}^\bot = \rad_{[n/2]+1}^\bot$.
\end{rem}

As in \S \ref{u11_E},
we define a $\sigma$-stable $\ri_F$-lattice $\J$ in $A$ by
(\ref{eq:def_J})
and an open compact subgroup $J$ 
of $G$ by $J = G\cap (1+\J)$.
For $b \in \g \cap \J^*$,
we denote by $\Psi_b$ the character of $J$
defined by (\ref{eq:Psi}).

\begin{lem}\label{lem:5c1_ad}
For $k \in \Z$,
the map
$\mathrm{ad}(\beta)$ induces an isomorphism
$\g_{k}^\bot/\g_{k+1}^\bot\simeq 
\g_{k-n}^\bot /\g_{k-n+1}^\bot$.
\end{lem}
\begin{proof}
Since $\{\seqrad_k\}_{k \in \Z}$ is  periodic,
it is sufficient to prove that the induced map is injective, 
for all $k \in \Z$.
For $X$ in $A$,
we use the notation
\begin{eqnarray*}
X 
= 
\left(
\begin{array}{cc}
X_{11} & X_{12}\\
X_{21} & X_{22}
\end{array}
\right),\
X_{ij} \in A^{ij}.
\end{eqnarray*}
Let $X \in \g_k^\bot$ satisfy $\ad(\beta)(X) \in \g_{k-n+1}^\bot$.
Then we obtain 
\[
\ad(\beta)(X) = 
\left(
\begin{array}{cc}
0 & -X_{12}\beta_2\\
\beta_2 X_{21} & \ad(\beta_2)(X_{22})
\end{array}
\right) \in \g_{k-n+1}^\bot
\]
and hence
$X_{12}\beta_2, \beta_2 X_{21}, \ad(\beta_2)(X_{22}) \in \rad_{k-n+1}$.
Since $\beta_2^{-1} \in \rad_{n}$,
we have 
$X_{12}$, $X_{21} \in \rad_{k+1}$
and
$X_{22} = X_{22} -s_2(X_{22})
= -\ad(\beta_2)(X_{22})\beta_2^{-1}/2 \in \rad_{k+1}$.
This completes the proof.
\end{proof}
We write $G' = G\cap B$, $J' = G'\cap J$,
$\Psi = \Psi_\beta$ and $\Psi' = \Psi|_J$.
Then, by Lemma~\ref{lem:5c1_ad},
we get the analogues of Propositions~\ref{prop:5c1_surj},
\ref{prop:5_2}, \ref{prop:E_surj}, and
Corollary~\ref{cor:5c1_J}.
\begin{thm}\label{thm:H-isom1}
With the notation as above,
there exists a support-preserving, $*$-isomorphism 
$\eta : \He(G'//J', \Psi') \simeq \He(G//J, \Psi)$.
\end{thm}
\begin{rem}
Since the centers of $G$ and $G'$ are compact,
the same remark as in the case (\ref{sec:h_int}c) holds,
that is,
$\eta$ induces a bijection from $\Irr(G')^{(J', 1)}$
to $\Irr(G)^{(J, \Psi)}$,
which preserves supercuspidality of representations.
\end{rem}

\begin{prop}\label{prop:H_inj}
For $g \in G'$, we have
$JgJ \cap G' = J'gJ'$.
\end{prop}
\begin{proof}
Put $M = (A^{11}\oplus A^{22})^\times$
and $x = 1_{V^1} - 1_{V^2} \in M$.
Since
$1 +\J$ is a pro-$p$ subgroup of $\widetilde{G}$
and stable under the adjoint action of $x$,
it follows from \cite{St2} Lemma 2.1
that
$(1+\J)g(1+\J)\cap M = ((1+\J)\cap M) g ((1 +\J)\cap M)$.

Consider the adjoint action of 
\[
h = 
\left(
\begin{array}{cc|cc}
0 & 0 & 0 & 1\\
0 &  1 &  0 & 0\\ \hline
0 & 0 & 1 & 0\\
\p d & 0 & 0 & 0
\end{array}
\right) \in B^\times
\]
on $M$.
Then $\Ad(h)$ induces an automorphism of 
$M$ of order 2
and $B^\times$ is the set of fixed points of this automorphism.
Note that $(1 + \J)\cap M$ is $Ad(h)$-stable.
Then, by \cite{St2} Lemma 2.1, we get
$(1+\J) g (1+\J)\cap B = ((1+\J)\cap B) g ((1 +\J)\cap B)$.

We note that
$B^\times$  is $\sigma$-stable and 
$(1+\J)\cap B$ is a $\sigma$-stable  pro-$p$ subgroup
of $B^\times$.
Since $G' = B \cap G$,
it follows from \cite{St2} Lemma 2.1 again that
$(1+\J) g (1+\J)\cap G' = J' g J'$.
This completes the proof.
\end{proof}

We put $B' = P_{\Lambda, 0}\cap G'$.
Then $B'$ is the Iwahori subgroup of $G'$ and
normalizes $J' = P_{\Lambda, n} \cap G'$.
We define elements $s_1$ and $s_2$ in $G'$ by
\[
s_1 = 
\left(
\begin{array}{cccc}
1 & 0 & 0 & 0\\
0 & 0 & 1 & 0\\
0 & 1 & 0 & 0\\
0 & 0 & 0 & 1
\end{array}
\right),\
s_2 = 
\left(
\begin{array}{cccc}
1 & 0 & 0 & 0\\ 
0 & 0 & \p^{-1} & 0\\
0 & \p & 0 & 0\\ 
0 & 0 & 0 & 1
\end{array}
\right).
\]
Put $S = \{s_1, s_2\}$ and $W' = \langle S \rangle$.
Then we have a Bruhat decomposition $G' = B' W' B'$.

\begin{lem}\label{lem:5c1_J'w}
Let $t$ be a non negative integer. 
Then

(i)
$[J' (s_1s_2)^t J': J'] = 
[J' (s_1s_2)^ts_1 J': J'] = 
[J'(s_2 s_1)^tJ' : J'] 
= q^{2t}$,
$[J' (s_2 s_1)^t s_2 J' : J'] = q^{2(t+1)}$.

(ii)
$[J (s_1s_2)^t J: J] = [J(s_2 s_1)^tJ : J]  = q^{6t},$
$[J (s_1s_2)^ts_1 J: J] = q^{6t+2}$,
$[J (s_2 s_1)^t s_2 J: J] = q^{6t+4}$.
\end{lem}
\begin{proof}
As in the proof of Lemma~\ref{lem:cosets},
we can directly calculate these indices using
(\ref{eq:filt-5c1}).
\end{proof}

For $\mu \in F$ and $\nu \in F^\times$,
 we set
\[
u(\mu) = 
\left(
\begin{array}{cccc}
1 & 0 & 0 & 0\\ 
0 & 1 & \mu & 0\\
0 & 0 & 1 & 0\\ 
0 & 0 & 0 & 1
\end{array}
\right),\
\underline{u}(\mu) = 
\left(
\begin{array}{cccc}
1 & 0 & 0 & 0\\ 
0 & 1 & 0 & 0\\
0 & \mu & 1 & 0\\ 
0 & 0 & 0 & 1
\end{array}
\right),\
h(\mu) = 
\left(
\begin{array}{cccc}
1 & 0 & 0 & 0\\ 
0 & \nu & 0 & 0\\
0 & 0 & \overline{\nu}^{-1} & 0\\ 
0 & 0 & 0 & 1
\end{array}
\right).
\]
For $g \in G'$,
let $e_g$ denote the element in 
$\He(G'//J', \Psi')$
such that
 $e_g(g) = 1$ and $\mathrm{supp}(e_g) = J'gJ'$.
\begin{thm}\label{thm:H_G'}
Suppose $m\geq 2$.
Then the algebra $\He(G'//J', \Psi')$ is generated by the elements $e_g$, $g \in B' \cup S$.
These elements are subject to the following relations:

(i) $e_k = \Psi(k) e_1$, $k \in J'$,

(ii) $e_k * e_{k'} = e_{kk'}$, $k, k' \in B'$,

(iii) $e_k * e_{s} = e_s * e_{sks}$, $s \in S$, $k \in B' \cap sB's$,

(iv) $e_{s_1}*e_{s_1} = e_1$,

$e_{s_2} * e_{s_2} = [J's_2 J' :J']\sum_{x \in \ri_0/\mi_0^2} 
e_{u(\p^{m-2}x\e)}$,

(v) For $\mu \in \ri_0^\times \e$,

$e_{s_1}* e_{u(\mu)}* e_{s_1}
= e_{u(\mu^{-1})} * e_{s_1} * e_{h(\mu)} * e_{u(\mu^{-1})}$,

$e_{s_2}* e_{\underline{u}(\p \mu)}* e_{s_2}
= [J's_2J':J']
e_{\underline{u}(\p \mu^{-1})} * e_{s_2} * e_{h(-\mu^{-1})} 
* e_{\underline{u}(\p \mu^{-1})}$.
\end{thm}
\begin{proof}
Note that $s_1$ normalizes $J'$.
Then 
the proof is very similar to that of
\cite{Harish} Chapter 3 Theorem 2.1.
\end{proof}

For $g \in G'$,
we denote by
 $f_g$ the element in $\He(G//J, \Psi)$ such that
$f_g(g) = 1$ and $\mathrm{supp}(f_g) = JgJ$.
\begin{thm}\label{thm:H_G}
Suppose $m \geq 2$.
Then the algebra
$\He(G//J, \Psi)$ is generated by $f_g$, $g \in B' \cup S$
and satisfies the following relations:

(i) $f_k = \Psi(k)f_1$, $k \in J'$,

(ii) $f_k * f_{k'} = f_{kk'}$, $k, k' \in B'$,

(iii) $f_k * f_{s} = f_s * f_{sks}$, $s \in S$, $k \in B' \cap sB's$,

(iv) $f_{s_1}*f_{s_1} = [Js_1J : J]f_1$,

$f_{s_2} * f_{s_2} = [Js_2 J :J]\sum_{x \in \ri_0/\mi_0^2} 
f_{u(\p^{m-2}x\e)}$,

(v) For $\mu \in \ri_0^\times \e$,

$f_{s_1}* f_{u(\mu)}* f_{s_1}
= -q f_{u(\mu^{-1})} * f_{s_1} * f_{h(\mu)} * f_{u(\mu^{-1})}$,

$f_{s_2}* f_{\underline{u}(\p \mu)}* f_{s_2}
=-q^3 f_{\underline{u}(\p \mu^{-1})} * f_{s_2} * f_{h(-\mu^{-1})} 
* f_{\underline{u}(\p \mu^{-1})}$.
\end{thm}
\begin{proof}
As in the proof of Theorem~\ref{thm:f},
it follows from Lemma~\ref{lem:5c1_J'w} that
$\He(G//J, \Psi)$ is generated by
$f_g$, $g \in S\cup B'$.
Relations (i), (ii), (iii) are obvious.
We shall prove relations (iv) and (v) on $s_2$ in the case
when $m = 2k+1$, $k \geq 1$.
The other cases can be treated in a similar fashion.

We will abbreviate $s = s_2$.
We can choose a common system of 
representatives for $J/J\cap sJs$
and $J\cap s Js\backslash J$ to be
\begin{eqnarray*} 
& x(a, A) = 
\left(
\begin{array}{cccc}
1 & 0 & 0 & 0\\
0 & 1 & 0 & 0\\
\p^{k+1}A & \p^m a\e & 1 & 0\\
0 & -\p^{k+1}\overline{A} & 0 & 1
\end{array}
\right),\ A \in \ri_F/\mi_F,\ a \in \ri_0/\mi_0^2.
\end{eqnarray*}
We note that $x(a,A)$ lies in the kernel of $\Psi$.

(iv)
As in the proof of Theorem~\ref{thm:f},
we obtain
\begin{eqnarray*}
f_s * f_s
 & = & 
 [JsJ:J] 
 \sum_{a\in \ri_0/\mi_0^2,\ A\in \ri_F/\mi_F} 
 f_1*\delta_{sx(a, A)s}*f_1.
\end{eqnarray*}
Since the $B$-component $sx(a,0)s$
of $sx(a,A)s$ lies in $\rad_0$,
it follows from
 Lemma~\ref{lem:intertwine_J}
 that
if $sx(a, A)s \in I_G(\Psi) = JG'J$,
then 
$sx(0, A)s \in 1 +\rad_{[n/2]+1}$ and hence $A \equiv 0$.
So $f_1*\delta_{sx(a, A)s}*f_1 = 0$
unless $A \equiv 0$.
Hence we have
\begin{eqnarray*}
f_s * f_s
 & = & 
 [JsJ:J] 
 \sum_{a\in \ri_0/\mi_0^2} 
 f_1*\delta_{sx(a, 0)s}*f_1\\
& = &
 [Js J :J]\sum_{a \in \ri_0/\mi_0^2} 
f_1*\delta_{u(\p^{m-2}a\e)}*f_1\\
& = &
[Js J :J]\sum_{a \in \ri_0/\mi_0^2} 
f_{u(\p^{m-2}a\e)}.
\end{eqnarray*}

(v)
Let $\mu \in \ri_0^\times \e$.
Put $u = \underline{u}(\p \mu) \in B'$.
As in the proof of Theorem~\ref{thm:f},
we have
\begin{eqnarray*}
f_s * f_u*f_s
 & = & 
 [JsJ:J] 
 \sum_{a\in \ri_0/\mi_0^2,\ A\in \ri_F/\mi_F} 
 f_1*\delta_{sx(a, A)us}*f_1.
\end{eqnarray*}
Put $x = x(0, A)$,
$\nu = \mu + \p^{m-1} a$,
$v = \underline{u}(\p \nu^{-1})$, and $h(-\nu^{-1})$.
Then
\begin{eqnarray*}
sx(a, A)us  
= sxs s\underline{u}(\p \nu) s
= sxs vshv
= [sxs, v]vshv (hv)^{-1}x(hv).
\end{eqnarray*}
Since $hv \in B'$,
we see that $(hv)^{-1}x(hv) \in J$ lies in the kernel of $\Psi$.
If we write $sxs = 1+y$ and $v = 1+z$,
then $y \in \rad_{8k-1}$ and $z \in \rad_6$.
So we have
$[sxs, v] \in P_{\Lambda, 8k+5} \subset P_{\Lambda, [n/2]+1}$
and
$[sxs, v] \equiv 1 + yzy \pmod{B^\bot +\rad_{n+1}}$.
These observations imply that
$[sxs, v] \in J$ and 
$\Psi([sxs, v]) = \Omega(\mathrm{tr}_{A/F_0}(\beta yzy))
= \Omega(2d\e \nu^{-1}A\overline{A})$.
We therefore have
\begin{eqnarray*}
f_s * f_u*f_s
& = &
 [JsJ:J] 
 \sum_{a\in \ri_0/\mi_0^2,\ A\in \ri_F/\mi_F} 
 \Omega^{-1}(2\e \nu^{-1}A\overline{A})
 f_1*\delta_{vshv}*f_1\\
 & = & 
-q  [JsJ:J] 
 \sum_{a\in \ri_0/\mi_0^2} 
 f_1*\delta_{vshv}*f_1\\
 & = & 
-q
 \sum_{a\in \ri_0/\mi_0^2} 
f _v*f_s*f_h*f_v \\
 & = & 
 -q^3  f_{\underline{u}(\p \mu^{-1})} * f_{s_2} * f_{h(-\mu^{-1})} 
* f_{\underline{u}(\p \mu^{-1})}.
\end{eqnarray*}
\end{proof}

\begin{rem}
If $m = 1$,
the algebra
$\He(G'//J', \Psi)$ is generated by the elements $e_g$, $g \in B' \cup S$.
In this case,
these elements are subject to the following relations:

(i) $e_k = \Psi(k)e_1$, $k \in J'$,

(ii) $e_k * e_{k'} = e_{kk'}$, $k, k' \in B'$,

(iii) $e_k * e_{s} = e_s * e_{sks}$, $s \in S$, $k \in B' \cap sB'S$,

(iv) $e_{s_1}*e_{s_1} = e_1$,

$e_{s_2} * e_{s_2} = 
(\sum_{y \in \ri_0\e/\mi_0\e,\ y \not \equiv 0}e_{s_2}
*e_{h(-y^{-1})} + q^2 e_1)(\sum_{x \in \ri_0/\mi_0} 
e_{u(x\e)})$,

(v) For $\mu \in \ri_0^\times \e$,

$e_{s_1}* e_{u(\mu)}* e_{s_1}
= e_{u(\mu^{-1})} * e_{s_1} * e_{h(\mu)} * e_{u(\mu^{-1})}$.

\noindent
We can easily see that the analogue of  Theorem~\ref{thm:H-isom1}
holds as well.
We omit the details.
\end{rem}

Combining 
Proposition~\ref{prop:H_inj},
the analog of Proposition~\ref{prop:E_surj}
and Theorems~\ref{thm:H_G'}, \ref{thm:H_G},
we complete the proof of Theorem~\ref{thm:H-isom1}.
Indeed, the map
\[
\eta(e_{s_1}) = -q^{-1} f_{s_1},\
\eta(e_{s_2}) = -q^{-1} f_{s_2},\
\eta(e_b) = f_b,\ b \in B'
\]
induces the required algebra isomorphism.

\subsection{Case (\ref{sec:h_int}a)}\label{sec:split}
Let 
$\Lambda$ denote the $\ri_F$-lattice sequence in 
(\ref{eq:st_20})
and $n = 2m-1$
a positive odd integer.
Let $\lambda \in \ri_F$ satisfy 
$\lambda \not\equiv \overline{\lambda} \pmod{\mi_F}$.
We consider the irreducible smooth representations of $G$
which contain
a skew stratum $[\Lambda, n, n-1, \beta]$,
where
\begin{eqnarray}
\beta = \p^{-m}
\left(
\begin{array}{cccc}
0 & 0 & 1 & 0\\
0 & 0& 0 &-1\\
\p \lambda & 0 & 0 & 0\\
0 & -\p \overline{\lambda} & 0 & 0
\end{array}
\right) \in \g_{\Lambda, -n}.
\end{eqnarray}
The assertions and proofs are very similar to those in \S 
\ref{u11_E} and \ref{5_1}.
So we shall be brief in this section.
\begin{rem}
By the argument in \S \ref{subsec:norm},
$\Irr(G)^{(P_{\Lambda, n}, \psi_\beta)}$
is precisely the set of equivalence classes of irreducible smooth representations of $G$
of level $n/2$
and of characteristic polynomial 
$(X -\lambda)^2(X-\overline{\lambda})^2 \pmod{\mi_F}$.
\end{rem}

Set $V^{-1} = F e_1 \oplus Fe_3$
and $V^1 = F e_2 \oplus F e_4$.
Then this stratum is split with respect to $V = V^{-1} \oplus V^1$.
For $j = \pm 1$,
we write 
$\beta_j = \beta|_{V^j}$, $E_j = F[\beta_j]$ and
$\Lambda^j(i) = \Lambda(i) \cap V^j$, $i \in \Z$.
Then Proposition~\ref{prop:pure} implies 
that $[\Lambda^j, n, n-1, \beta_j]$, $j = \pm 1$,
are simple strata.
We set $A^{ij} = \mathrm{Hom}_F(V^j, V^i)$, $i, j = \pm 1$.
Since $\beta_j^2 \in F^\times \cdot 1_{V^j}$,
we obtain an $(E_j, E_j)$-bimodule projection
$s_j : A^{jj} \rightarrow E_j$ by
\[
s_j(X) = (X + \beta_j X \beta_j^{-1})/2,\
X \in A^{jj}.
\]
Writing $B$ for the $A$-centralizer of $\beta$ 
and $B^\bot$ for its orthogonal component with respect to
the pairing induced by $\mathrm{tr}_{A/F}$,
we have
\[
B = E_{-1} \oplus E_1,\
B^\bot = \ker s_{-1} \oplus \ker s_1 \oplus A^{-1, 1} 
\oplus A^{1, -1},
\]
and $A = B \oplus B^\bot$.
For $j = \pm 1$,
the map $s_j$ satisfies $s_j(\rad_k(\Lambda^j)) = 
\rad_k(\Lambda^j)\cap E_j$.
Combining this with \cite{BK2} (2.9),
we obtain 
\begin{eqnarray*}
\rad_k(\Lambda) = \rad_k(\Lambda)\cap B
\oplus \rad_k(\Lambda)\cap B^\bot,\ 
k \in \Z.
\end{eqnarray*}
So
we can define a $\sigma$-stable $\ri_F$-lattice $\J$
in $A$ by
(\ref{eq:def_J})
and an open compact subgroup $J = \J\cap G$ of $G$.
For $b \in \g \cap \J^*$,
we denote by $\Psi_b$ the character of $J$
defined by (\ref{eq:Psi}).

As in \S \ref{u11_E},
we abbreviate $\rad_k = \rad_k(\Lambda)$,
$\rad_k' = \rad_k \cap B$, $\rad_k^\bot = \rad_k \cap B^\bot$, etc.
\begin{lem}\label{lem:split_ad}
For  $k \in \Z$,
$\ad(\beta)$ induces an isomorphism
$\g_{k}^\bot/\g_{k+1}^\bot \simeq 
\g_{k-n}^\bot /\g_{k-n+1}^\bot$.
\end{lem}
\begin{proof}
It follows from \cite{BK2} (3.7) Lemma 1 that
$\ad(\beta)$ maps $\rad_k \cap A^{-1,1}$ onto
$\rad_{k-n}\cap A^{-1, 1}$
and $\rad_k \cap A^{1,-1}$ onto
$\rad_{k-n}\cap A^{1, -1}$.

Let $j \in \{ \pm 1\}$.
For $x \in \rad_{k-n} \cap \ker s_j$,
we have $-x\beta_j^{-1}/2 \in \rad_{k} \cap \ker s_j$
and 
$\ad(\beta_j)(-x\beta_j^{-1}/2)
= (x- \beta_j x\beta_j^{-1} )/2
= x -s_j(x) = x$.
This implies
 that $\ad(\beta_j)$ maps $\rad_k \cap \ker s_j$ onto
$\rad_{k-n}\cap \ker s_j$,
so that
$\ad(\beta)$ maps $\rad_k^\bot$ onto $\rad_{k-n}^\bot$.

By the periodicity of $\{\rad_k \}_{k \in \Z}$,
we conclude that 
$\ad(\beta)$ induces an isomorphism
$\rad_{k}^\bot/\rad_{k+1}^\bot \simeq 
\rad_{k-n}^\bot /\rad_{k-n+1}^\bot$.
The lemma follows immediately from this.
\end{proof}

We also set $G' = G\cap B$, $J' = G'\cap J$,
$\Psi = \Psi_\beta$, and $\Psi' = \Psi|_J$.
Then, by Lemma~\ref{lem:split_ad},
we get the analogues of Propositions~\ref{prop:5_2} and  \ref{prop:E_surj}, that is,
$\Irr(G)^{(P_{\Lambda, n}, \psi_\beta)}
= \Irr(G)^{(J, \Psi)}$
and $I_G(\Psi) = JG'J$.
Applying the proof of 
Proposition~\ref{prop:H_inj}, 
we have
$JgJ \cap G' = J'gJ'$, for $g \in G'$.

We shall establish a Hecke algebra isomorphism
$\eta : \He(G'//J', \Psi') \simeq \He(G//J, \Psi)$.
Since $G'$ is abelian,
the structure of the algebra $\He(G'//J', \Psi')$ is obvious.

We have $G' = \langle \zeta \rangle B'$,
where $B' = G'\cap P_{\Lambda, 0}$ and
\begin{eqnarray}
\zeta = 
\left(
\begin{array}{cccc}
0 & 0 & \p^{-1}\lambda^{-1} & 0\\
0 & 0 & 0 & 1\\
1 & 0 & 0 & 0\\
0 & \p \overline{\lambda} & 0 & 0
\end{array}
\right).
\end{eqnarray}
The group $B'$ normalizes the pair $(J, \Psi)$
and $\zeta$ satisfies the following:
\begin{lem}\label{lem:split_cosets}
For $t \in \Z$,
we have
$[J  \zeta^t J: J] = q^{4|t|}$.
\end{lem}
\begin{proof}
It follows from direct computation.
\end{proof}

For $g \in G'$,
let $e_g$ denote the element in 
$\He(G'//J', \Psi')$
such that
 $e_g(g) = 1$ and $\mathrm{supp}(e_g) = J'gJ'$,
and
let $f_g$ denote the element in $\He(G//J, \Psi)$ such that
$f_g(g) = 1$ and $\mathrm{supp}(f_g) = JgJ$.
\begin{lem}\label{lem:split_relations}
$f_{\zeta^{-1}}*f_\zeta =
f_\zeta * f_{\zeta^{-1}} = q^4 f_1$.
\end{lem}
\begin{proof}
As in the proof of Theorem~\ref{thm:f},
we have
\begin{eqnarray*}
f_{\zeta^{-1}}* f_\zeta
& = & [J\zeta J:J]
\sum_{y \in J/J\cap \zeta J\zeta^{-1}}
\Psi^{-1}(y) f_1 *\delta_{\zeta^{-1}y \zeta}*f_1.
\end{eqnarray*}
Suppose that $m =2k$.
Then we can choose a system of representative
for $J/J\cap \zeta J\zeta^{-1}$ to be
\begin{eqnarray*}
x(a, b, A)= 
\left(
\begin{array}{cccc}
1 & 0 & 0 & 0\\
\p^k A & 1 & \p^k a\e & 0\\
0 & 0 & 1 & 0\\
\p^{k+1}b\e & 0 & -\p^k \overline{A} & 1
\end{array}
\right),\
a, b \in \ri_0/\mi_0,\ A \in \ri_F/\mi_F.
\end{eqnarray*}
The $B$-component of $\zeta^{-1}x(a, b, A)\zeta$ lies in 
$\rad_0$,
so that, if $\zeta^{-1}x(a, b, A)\zeta \in JG'J$,
then by  Lemma~\ref{lem:intertwine_J},
the $B^\bot$-part of $\zeta^{-1}x(a, b, A)\zeta$
belongs to $\rad_{[n/2]+1}$.
We conclude that
$f_1 *\delta_{\zeta^{-1}x(a, b, A) \zeta}*f_1 = 0$ unless $a \equiv b \equiv A \equiv 0$.
So we have
$f_{\zeta^{-1}}* f_\zeta
 = q^4 f_1$, as required.
 
The case of $m = 2k+1$, $k \geq 0$ 
and the relation 
$f_\zeta * f_{\zeta^{-1}} = q^4 f_1$
are quite similar to this. 
\end{proof}

\begin{thm}\label{thm:slpit_isom}
With the notation as above,
there is a  $*$-isomorphism 
$\eta : \He(G'//J', \Psi') \simeq \He(G//J, \Psi)$
with support preservation.
\end{thm}
\begin{proof}
Since $I_G(\Psi) = JG' J$,
the algebra
$\He(G//J, \Psi)$ is spanned by $f_g$, $g \in G'$.
Since $G' = \langle \zeta \rangle B'$,
we can write $g = \zeta^t b$ where $t \in \Z$, $b \in B'$.
By Lemma~\ref{lem:split_cosets},
we have
$f_g = f_\zeta^t *f_b$ if $t \geq 0$
and $f_g = f_{\zeta^{-1}}^{|t|}* f_b$ otherwise.

Since $B'$ normalizes the pair $(J, \Psi)$,
we have $f_b * f_{g} = f_{bg} = 
f_{gb} = f_{g} * f_b$, for $b\in B'$ and $g \in G'$.
In particular,
$f_b$ lies in the center of $\He(G//J, \Psi)$,
for $b \in B'$.
Therefore the map
\[
\eta(e_\zeta) = q^{-2}f_\zeta,\
\eta(e_{\zeta^{-1}}) = q^{-2}f_{\zeta^{-1}},\
\eta(e_b) = f_b,\ b \in B',
\]
induced the required isomorphism
$\eta : \He(G'//J', \Psi') \simeq \He(G//J, \Psi)$
by Lemma~\ref{lem:split_relations}.
\end{proof}

\begin{rem}
As in Remark~\ref{rem:u11_E},
the map $\eta$ induces a bijection 
from $\Irr(G')^{(J', 1)}$
to $\Irr(G)^{(J, \Psi)}$.
Since $G'$  is abelian and not compact,
every irreducible smooth representation of $G$
which contains $(J, \Psi)$ is not supercuspidal.
\end{rem}

\end{document}